\documentclass[a4paper,11pt]{article}

	\usepackage{titling}
	\pretitle{\begin{center}\Large\bfseries} 
	\posttitle{\par\end{center}\vskip 0.5em}

	\usepackage{amsfonts,amssymb,amsmath,amsthm}
	\usepackage{graphicx}
	\graphicspath{{figures/}}
	\usepackage[font=small, labelfont={bf, sf}]{caption}
	
	\usepackage{subcaption}
	\usepackage{mathtools}
	\usepackage[dvipsnames]{xcolor}
	\usepackage[inline, shortlabels]{enumitem}
	\usepackage[T1]{fontenc}
	\usepackage[utf8]{inputenc}
	\usepackage{bm}
	\usepackage{xparse}
	\usepackage{dsfont}
	\usepackage{mathrsfs}
	\usepackage{float}
	
	\usepackage{authblk}
	\usepackage{indentfirst}
	\usepackage[sf, bf, medium]{titlesec}

	\usepackage[margin=2.5cm]{geometry}

	\usepackage{tikz}
	\usetikzlibrary{
		arrows, 
		arrows.meta, 
		patterns, 
		calc, 
		positioning,
		decorations.pathreplacing,
		intersections,
		shapes
		}

	\usepackage[
	    backend=bibtex,
	    style=numeric,
	    	sortcites=true,
	    sorting=nty,
	    giveninits=true,
	    maxbibnames=99
		]{biblatex}
	\DeclareFieldFormat{labelnumberwidth}{#1\adddot}
	
	\renewbibmacro*{in:}{}
	\DeclareFieldFormat[article]{title}{#1}
	\DeclareFieldFormat[article]{volume}{\mkbibbold{#1}}

	\usepackage[bookmarks]{hyperref}
	\hypersetup{
	    colorlinks=true,
		citecolor=[rgb]{0,0.173,0.620},
		linkcolor=[rgb]{0,0.173,0.620},
		urlcolor=[rgb]{0,0.173,0.620},
		pdfpagemode=UseNone,
		pdfstartview=FitH,
		pdfdisplaydoctitle=true,
		pdftitle={Sharp global L2 estimates for oscillatory integral operators with quasi-homogeneous phases},
		pdfauthor={Minbeom Kang, Joonil Kim},
		pdflang=en-US
	 }
	\usepackage{cleveref}

	\theoremstyle{plain}
	\newtheorem{theorem}{\sffamily Theorem}\numberwithin{theorem}{section}
	\newtheorem{lemma}{\sffamily Lemma}\numberwithin{lemma}{section}
	\newtheorem{proposition}{\sffamily Proposition}\numberwithin{proposition}{section}
	\newtheorem{corollary}{\sffamily Corollary}\numberwithin{corollary}{section}

	\theoremstyle{definition}
	\newtheorem{definition}{\sffamily Definition}\numberwithin{definition}{section}
	
	\newtheoremstyle{example}
	  {\topsep}
	  {\topsep}
	  {\normalfont}
	  {}
	  {\bfseries \sffamily \itshape}
	  {}
	  {0.5em}
	  {}
	 
	\theoremstyle{example}
	\newtheorem{example}{Example}\numberwithin{example}{section}
	\newtheorem*{example*}{Examples~\ref{ex:1@intro}--\ref{ex:3@intro} revisited.}
	
	\theoremstyle{remark}
	\newtheorem{remark}{\sffamily Remark}\numberwithin{remark}{section}
	\newtheorem*{ack*}{Acknowledgment}

	\numberwithin{equation}{section}	

	\newtheoremstyle{mainthmstyle}
	  {\topsep}   
	  {\topsep}   
	  {\itshape}  
	  {}          
	  {\bfseries \sffamily} 
	  {.}         
	  {.5em}     
	  {#1\thmnote{ (#3)}} 

	\theoremstyle{mainthmstyle}
	\newtheorem{main}{Main Theorem}
	\crefname{main}{Main Theorem}{Main Theorems}
	\Crefname{main}{Main Theorem}{Main Theorems}
	\crefformat{main}{#2\textcolor{black}{Main Theorem}#3}
	\Crefformat{main}{#2\textcolor{black}{Main Theorem}#3}

\newcommand{\C}{\mathbb{C}}

\newcommand{\N}{\mathbb{N}}

\renewcommand{\P}{\mathbb{P}}
\newcommand{\Q}{\mathbb{Q}}
\newcommand{\R}{\mathbb{R}}

\newcommand{\Z}{\mathbb{Z}}

\newcommand\fN{{\mathfrak {N}}}

\newcommand\fb{{\mathfrak {b}}}

\newcommand\fm{{\mathfrak {m}}}

\newcommand\fp{{\mathfrak {p}}}
\newcommand\fq{{\mathfrak {q}}}

\newcommand\cC{{\mathcal {C}}}

\newcommand\cL{{\mathcal {L}}}

\newcommand\cR{{\mathcal {R}}}
\newcommand\cS{{\mathcal {S}}}
\newcommand\cT{{\mathcal {T}}}

\newcommand\1{{\bf 1}}

\newcommand\bC{{\bf {C}}}

\newcommand\bN{{\bf {N}}}

\newcommand\be{{\bf {e}}}

\newcommand\bu{{\bf {u}}}
\newcommand\bv{{\bf {v}}}

\newcommand\by{{\bf {y}}}

\newcommand\tJ{{\widetilde {J}}}
\newcommand\tj{{\widetilde j}}

\newcommand{\toinf}{\to\infty}

\DeclareMathOperator{\conv}{\mathbf{conv}}
\DeclareMathOperator{\cone}{cone}
\DeclareMathOperator{\supp}{supp}
\DeclareMathOperator{\cl}{cl}
\DeclareMathOperator*{\esssup}{ess\,sup}

\def\plane#1#2{\pi_{#1, #2}}

\newcommand\df{\delta_\mathrm{ for}}
\newcommand\db{\delta_\mathrm{ bac}}
\newcommand\loc{\mathrm{ loc}}
\newcommand\glo{\mathrm{ glo}}
\newcommand\op{L^2 \to L^2} 
\renewcommand{\Re}{\operatorname{Re}}

\def\dx{\mathop{}\!\mathrm{d}x}
\def\dy{\mathop{}\!\mathrm{d}y}
\def\d#1{\mathop{}\!\mathrm{d}#1}

\newcommand{\defeq}{\coloneqq}
\newcommand{\eqdef}{\eqqcolon}

\newcommand{\Def}[2]{%
	\begin{definition}[#1]
	{#2}
	\end{definition}%
	}

\newcommand{\Thm}[2]{%
	\begin{theorem}[#1]
	{#2}
	\end{theorem}%
	}

\newcommand{\Ex}[2]{%
	\begin{example}[#1]
	{#2}
	\end{example}%
	}

\newcommand{\Lemma}[2]{%
	\begin{lemma}[#1]
	{#2}
	\end{lemma}%
	}

\newcommand{\Prop}[2]{%
	\begin{proposition}[#1]
	{#2}	
	\end{proposition}%
	}

\newcommand{\pf}[1]{%
	\begin{proof}
	{#1}	
	\end{proof}%
	}

\newcommand{\Pf}[2]{%
	\begin{proof}[#1]
	{#2}	
	\end{proof}%
	}

\newcommand{\Cor}[2]{%
	\begin{corollary}[#1]
	{#2}	
	\end{corollary}%
	}

\newcommand{\Rmk}[2]{%
	\begin{remark}[#1]
	{#2}	
	\end{remark}%
	}

\DeclarePairedDelimiter{\paren}{(}{)}
\DeclarePairedDelimiter{\bracket}{[}{]}
\DeclarePairedDelimiter{\set}{\{}{\}}

\DeclarePairedDelimiter{\abs}{\lvert}{\rvert}
\DeclarePairedDelimiter{\norm}{\lVert}{\rVert}
\DeclarePairedDelimiter{\inn}{\langle}{\rangle}

\newcommand{\ii}{\mathrm{i}}
\newcommand{\ee}{\mathrm{e}}

\title{Sharp global \(L^2\) estimates for oscillatory integral operators with quasi-homogeneous phases} 

\author[1]{Minbeom Kang}
\author[1]{Joonil Kim}
\affil[1]{Department of Mathematics, Yonsei University, Seoul 03722, Korea}
\date{}

\begin{document}

\allowdisplaybreaks

\maketitle

\begin{abstract}
For one-dimensional oscillatory integral operators with real-analytic phases, the seminal work of Phong and Stein established that the Newton polygon of the mixed Hessian of the phase determines the sharp local \(L^2\)-decay. The local theory relies on spatial localization near the origin and therefore does not capture the global geometry of the real zero set of the Hessian. For globally defined operators, unbounded real zero curves of the Hessian can make the global decay strictly slower than the local decay or destroy \(L^2\)-boundedness altogether.

We establish a global analogue of the Phong--Stein theory for quasi-homogeneous polynomial phases. We introduce a family of Newton polygons adapted to the regions near and away from the real zero curves of the Hessian, together with the corresponding box-area scales. The relative geometry of these Newton polygons and area scales characterizes global \(L^2\)-boundedness and determines the sharp \(L^2\)-estimates in the high- and low-frequency regimes. As an application, we determine the exact range of anisotropic homogeneous \(L^2\)-Sobolev estimates for the associated global Radon transform.
\end{abstract}

\bigskip

\noindent
\textbf{Keywords } Oscillatory integral operators \(\cdot\) Newton polygons \(\cdot\) global \(L^2\)-estimates

\medskip
\noindent
\textbf{Mathematics Subject Classification } Primary 42B20 \(\cdot\) 47G10

\bigskip


\section{Introduction} \label{sec:intro}

In this paper, we study the globally defined oscillatory integral operator 
	\begin{equation} \label{eq:T_glo}
		\cT^\lambda_{\glo} f(x)\defeq \int \ee^{\ii\lambda P(x,y)} f(y)\dy 
		\quad \text{for } f \in \cS(\R),
	\end{equation}
where \(P \in \R[x,y]\) is a polynomial phase and \(\lambda>0\). Throughout this paper, we assume \(P(x,0)=P(0,y)=0\), as pure \(x\)- and \(y\)-terms do not affect the \(L^2\)-operator norm. We aim to characterize global \(L^2\)-boundedness and determine the sharp \(L^2\)-estimates. We fully resolve this problem for quasi-homogeneous polynomial phases.

We first recall the local theorem of Phong and Stein, which provides the starting point for our study. Let \(S(x,y)\) be a real-analytic phase function whose Hessian \(S_{xy}\) does not vanish identically. Phong and Stein considered the operator \(\cT^\lambda_\loc\) defined as
	\begin{equation} \label{eq:T_loc}
		\cT^\lambda_{\loc}f(x)=\int \ee^{\ii\lambda S(x,y)}\chi(x,y)f(y)\dy,
	\end{equation}
where \(\chi\in \cC_0^\infty(\R^2)\) is a smooth cutoff function supported in a sufficiently small neighborhood of the origin. Let \(S_{xy}(x,y)=\sum_{p,q} c_{p,q}x^py^q\) be the Taylor expansion of \(S_{xy}\) at the origin, and let \(\Lambda(S_{xy})\defeq \set{(p,q) \in \Z_{\ge0}^2 : c_{p,q}\neq 0}\) be its Taylor support.
Writing \(\be_1=(1,0)\) and \(\be_2=(0,1)\), we express the local model domain \([-1,1]^2\) as
	\[
		D_{\set{\be_1, \be_2}}
		= \cl\set[\big]{(x,y)\in (\R\setminus\set{0})^2: |x|^1|y|^0\le 1,\ |x|^0|y|^1\le 1},
	\]
where \(\cl\) denotes closure in \(\R^2\). In this notation, the \emph{local Newton polygon} \(\bN(S_{xy})\) takes the form
	\begin{equation} \label{eq:local_NP}
		\bN(S_{xy})
		= \conv\paren[\big]{\Lambda(S_{xy})+\R_{\ge0}^2}
		= \conv\paren[\big]{\Lambda(S_{xy})+\cone\set{\be_1,\be_2}}.
	\end{equation}
We use \(\conv E\) to denote the convex hull of a set \(E\), and set \(\cone B \defeq \set[\big]{\sum_{\fb \in B} \alpha_\fb \fb : \alpha_\fb \ge 0}\).
The \emph{local Newton distance} \(\delta_{\loc}\) is the number satisfying
	\begin{equation} \label{eq:local_newton_distance}
		\bN(S_{xy})\cap \cone(\1)
		= [\delta_{\loc},\infty)\1, \quad \text{where } \1=(1,1).
	\end{equation}
Here and below, for an interval \(I \subset \R\) and \(\bu \in \R^2\), we write \(I\bu \defeq \set{t\bu: t \in I}\). The theorem of Phong and Stein then reads as follows.
\Thm{Phong--Stein \cite{PS97}}{\label{thm:PS@intro}
Let \(\cT^\lambda_\loc\) be defined as in \eqref{eq:T_loc}. If the support of \(\chi\) is sufficiently small, then
	\[
		\norm{\cT^\lambda_\loc}_{L^2(\R)\to L^2(\R)}
		\lesssim \lambda^{-\frac{1}{2(\delta_{\loc}+1)}}
		\quad \text{for } \lambda \ge 1.
	\]
Furthermore, if \(\chi(0,0)\neq0\), this estimate is sharp.
}

A local dyadic block with \(|x| \approx 2^{-j_1}\) and \(|y| \approx 2^{-j_2}\) has area \(\approx (2^{-j_1})^1(2^{-j_2})^1\). The vector \(\1=(1,1)\) in \eqref{eq:local_newton_distance} is the exponent vector of this area scale, thereby determining the diagonal ray \(\cone(\1)\) along which the local Newton distance is measured.

\subsection{The global problem} \label{subsec:glo_prob}

We now compare the globally defined operator in \eqref{eq:T_glo} with its local counterpart. For the local comparisons below, we take \(\chi(0,0) \neq 0\) in \eqref{eq:T_loc}. We begin with the simplest global example, provided by the Fourier transform. For \(P(x,y)=xy\), Plancherel's theorem gives \(\norm{\cT^\lambda_\glo}_{\op} \approx \lambda^{-1/2}\), which agrees with the local estimate corresponding to \(\delta_\loc =0\).

We present the following examples to demonstrate three global phenomena: 
\begin{enumerate}[label=\rm (\arabic*)]
	\item 
		local \(L^2\)-boundedness does not imply global \(L^2\)-boundedness;
	\item
		the local Newton distance does not determine the global \(L^2\)-decay;
	\item
		the sharp global \(L^2\)-estimates differ between the regimes \(\lambda  \ge 1\) and \(0< \lambda < 1\).
\end{enumerate}

\Ex{}{\label{ex:1@intro}
Consider \(P(x,y)=-x^2y\). Since \(P_{xy}(x,y)=-2x\), in view of \labelcref{eq:local_NP,eq:local_newton_distance}, one has
	\[
		\bN(P_{xy}) 
		= \conv\paren[\big]{\set{(1,0)} + \cone\set{\be_1,\be_2}} \quad \text{and} \quad
		\bN(P_{xy}) \cap \cone(\1) = [1,\infty)\1,
	\]
yielding \(\delta_\loc =1\). Hence, \Cref{thm:PS@intro} gives the sharp local estimate
	\[
		\norm{\cT^\lambda_\loc}_{\op} \approx \lambda^{-\frac{1}{4}} 
		\quad \text{for } \lambda \ge1.
	\]
Globally, the operator acts as \(\cT^\lambda_\glo f(x)= \widehat f(\lambda x^2)\). Plancherel's theorem and the change of variables \(u=\lambda x^2\) yield
	\[
		\norm{\cT^\lambda_\glo f}_2^2 
		\approx \lambda^{-\frac{1}{2}}\int_0^\infty |\widehat f(u)|^2 u^{-1/2}\d{u}
		\quad \text{for } f \in \cS(\R).
	\]
For each \(\epsilon>0\), choose a nonzero \(f_\epsilon \in \cS(\R)\) such that \(\widehat f_\epsilon \in \cC_0^\infty((0, \epsilon))\). Then
	\[
		\frac{\norm{\cT^\lambda_\glo f_\epsilon}_2^2}
		{\norm{f_\epsilon}_2^2}
		\gtrsim \lambda^{-1/2}\epsilon^{-1/2}
		\longrightarrow +\infty \quad\text{as } \epsilon\to 0.
	\]
Thus \(\cT^\lambda_\glo\) has no bounded extension to \(L^2(\R)\), although its localized counterpart is bounded. Therefore, the local Newton polygon alone does not characterize global \(L^2\)-boundedness of \(\cT^\lambda_\glo\).
}

\Ex{}{\label{ex:2@intro}
Let \(P(x,y)=x^2y^2-4xy\). Since \(P_{xy}(x,y)=4xy-4\), in view of \labelcref{eq:local_NP,eq:local_newton_distance}, one has
	\[
		\bN(P_{xy}) 
		= \conv\paren[\big]{\set{(0,0),(1,1)}+\cone\set{\be_1,\be_2}} \quad \text{and} \quad
		\bN(P_{xy}) \cap \cone(\1) = [0,\infty)\1,
	\]
yielding \(\delta_\loc =0\). Hence \Cref{thm:PS@intro} gives the sharp local estimate
	\[
		\norm{\cT^\lambda_\loc}_{\op} \approx \lambda^{-\frac{1}{2}}
		\quad \text{for } \lambda \ge1.
	\]

Although \(P_{xy}\) is nondegenerate at the origin, it vanishes along the curve \(xy=1\). 
A Knapp-type argument directly gives a lower bound. Consider an axis-parallel box centered at \((1,1)\) with side lengths comparable to \(\lambda^{-1/3}\). Testing \(\cT^\lambda_\glo\) against suitably modulated functions supported on its coordinate projections gives
	\[
		\norm{\cT^\lambda_\glo}_{\op} 
		\gtrsim \lambda^{-\frac{1}{3}} \quad \text{for } \lambda \ge 1.
	\]
This lower bound already rules out the local rate \(\lambda^{-1/2}\). The loss of decay comes from the unbounded zero curve \(xy=1\), which is not encoded by the local Newton polygon. Therefore, the local Newton distance does not determine the sharp global high-frequency decay of \(\cT^\lambda_\glo\).
}

\Ex{}{\label{ex:3@intro}
To illustrate the third phenomenon, we revisit the phase \(P(x,y)=x^2y^2-4xy\) from \Cref{ex:2@intro}. The same Knapp-type argument, now applied to a box centered at \((\lambda^{-1/4}, \lambda^{-1/4})\) with side lengths comparable to \(\lambda^{-1/4}\), gives
	\[
		\norm{\cT^\lambda_\glo}_{\op} 
		\gtrsim \lambda^{-\frac{1}{4}} \quad \text{for } 0< \lambda < 1.
	\]
	
As \Cref{thm:main} and \Cref{cor:full_operator} will show, the lower bounds in \Cref{ex:2@intro,ex:3@intro} are both sharp. Consequently,
	\[
		\norm{\cT^\lambda_\glo}_{\op}
		\approx
		\begin{cases}
			\lambda^{-1/3}, & \text{if } \lambda\ge1,\\[0.5ex]
			\lambda^{-1/4}, & \text{if } 0<\lambda<1.
		\end{cases}
	\]
Thus the sharp \(L^2\)-estimates have different exponents of \(\lambda\) in the high- and low-frequency regimes.
}

Together, these examples show that the cutoff-free problem requires a global geometric framework that determines both \(L^2\)-boundedness and the sharp \(L^2\)-estimates in the high- and low-frequency regimes. This observation leads to the following questions:
\begin{enumerate}[wide=0pt, label= \bfseries \sffamily Question \arabic*.]
	\item For a fixed \(\lambda>0\), which polynomial phases \(P\) yield a bounded operator \(\cT^\lambda_\glo\)  on \(L^2(\R)\), and which do not? Which Newton polygons determine this dichotomy between boundedness and unboundedness?
	\item When \(\cT^\lambda_\glo\) is globally bounded on \(L^2(\R)\), what are its sharp \(L^2\)-estimates in the high- and low-frequency regimes? Which Newton polygons determine the estimate in each regime?
\end{enumerate}

\subsection{Main results} \label{subsec:main_results}

We answer these questions for quasi-homogeneous polynomial phases \(P\) with \(P_{xy}\not\equiv 0\). A nonzero polynomial is \emph{quasi-homogeneous} if its Taylor support is contained in an affine line. We allow this line to have arbitrary direction. As \Cref{subsec:scaling} shows, the direction of this line determines whether the global problem can be reduced to the local theory by scaling. We develop a global framework for all such phases, including those for which this reduction fails.

We first introduce Newton polygons adapted to regions defined by finitely many monomial inequalities.\footnote{Regions of this type were studied by Nagel and Pramanik \cite{NP09} in connection with maximal averages.} Viewing \(\bN(P_{xy}) = \conv\paren[\big]{\Lambda(P_{xy})+\cone\set{(1,0),(0,1)}}\) as the Newton polygon adapted to the local region \([-1,1]^2\), we extend the same construction to the adapted regions defined below, following the domain-adapted Newton polyhedra of \cite{Kim23}.

\Def{Cones and adapted regions}{\label{def:cone_domain}
Let \(B \subset \Q^2\) be a finite nonempty set.
\begin{enumerate}[(i)]
	\item 
		The cone generated by \(B\) is \(\cone B = \set[\big]{\sum_{\fb \in B} \alpha_\fb \fb : \alpha_\fb \ge 0}\).
	\item
		For a collection \(\bC = (C_\fb)_{\fb \in B}\) of positive constants, define
			\begin{equation} \label{eq:domain}
				D_{B,\bC} \defeq 
				\cl \set[\big]{(x,y)\in (\R\setminus\set{0})^2: |x|^{b_1}|y|^{b_2}\le C_\fb \text{ for all } \fb =(b_1,b_2)\in B}.
			\end{equation}
		When \(C_\fb = 1\) for every \(\fb \in B\), we simply write \(D_B\).
	\item
		For a nonempty region \(D\subset \R^2\), we write \(D \simeq D_B\) if there exist collections \(\bC=(C_\fb)_{\fb \in B}\) and \(\bC^\prime=(C^\prime_\fb)_{\fb \in B}\) of positive constants such that \(D_{B,\bC^\prime} \subset D \subset D_{B,\bC}\).
\end{enumerate}
Here and throughout this paper, we consider only nonempty regions \(D_{B,\bC}\).
}

\Def{Domain-adapted Newton polygon}{\label{def:NP}
Let \(P \in \R[x,y]\) be a nonzero real polynomial with Taylor support \(\Lambda(P)\), and let \(D \simeq D_B\). The \emph{Newton polygon of \(P\) adapted to \(D\)} is
	\begin{equation} \label{eq:NP}
		\bN(P,D)
		\defeq \conv\paren[\big]{\Lambda(P)+\cone B},
	\end{equation}
where \(\conv E\) denotes the convex hull of a set \(E \subset \R^2\). For a finite fractional power series \(F(x,y)\), we define \(\Lambda(F)\) and \(\bN(F,D)\) in the same manner.
}

\Rmk{}{
The Newton polygon \(\bN(P,D)\) is well defined. Namely, if \(D \simeq D_B\) and \(D \simeq D_{B^\prime}\), then \(\cone B = \cone B^\prime\). To see this, it suffices to verify that if \(D_{B^\prime, \bC^\prime} \subset D_{B, \bC}\), then \(\cone B \subset \cone B^\prime\). By duality for finitely generated cones, this is equivalent to
	\[
		\set{\fq \in \R^2 : \fq \cdot \fb^\prime \ge 0 \text{ for all } \fb^\prime \in B^\prime}
		\subset \set{\fq \in \R^2 : \fq \cdot \fb \ge 0 \text{ for all } \fb \in B}.
	\]
Fix \(\fq=(q_1,q_2)\) in the set on the left and \((x_0,y_0) \in D_{B^\prime, \bC^\prime}\) with \(x_0y_0 \neq 0\). Then, for all \(\alpha \ge 0\), one has \((2^{-\alpha q_1}x_0, 2^{-\alpha q_2}y_0) \in D_{B^\prime, \bC^\prime} \subset D_{B,\bC}\), and hence \(2^{-\alpha(\fq\cdot \fb)}|x_0|^{b_1}|y_0|^{b_2} \le C_\fb\) for all \(\fb = (b_1, b_2) \in B\). Letting \(\alpha\toinf\) yields \(\fq\cdot \fb \ge0\). Consequently, \(\bN(P,D)\) does not depend on \(B\).
}

\Ex{}{\label{ex:NP}
For a nonzero polynomial \(P\), the Newton polygons adapted to the model domains \([-1,1]^2=D_{\set{\be_1, \be_2}}\) and \(\R^2=D_{\set{\bm 0}}\) are given as follows.
\begin{itemize}
	\item \(\bN(P,D_{\set{\be_1, \be_2}})=\conv\paren[\big]{\Lambda(P)+\cone\set{\be_1, \be_2}}\) recovers the local Newton polygon in \eqref{eq:local_NP}.
	\item \(\bN(P,D_{\set{\bm 0}})=\conv\paren[\big]{\Lambda(P)}\) is the convex hull of the Taylor support.
\end{itemize}
In particular, if \(D=\cl\set[\big]{(x,y) \in (\R\setminus\set{0})^2:|x|^{-3}|y| \le 4}\), then \(D \simeq D_{\set{(-3,1)}}\), and therefore 
	\[
		\bN(P, D) = \conv\paren[\big]{\Lambda(P)+\cone\set{(-3,1)}}.
	\]
}

With this geometric framework in hand, we turn to the global operator \(\cT^\lambda_\glo\) in \eqref{eq:T_glo}. We decompose it spatially into four quadrants:
	\begin{equation} \label{eq:T_quadrant_decom}
	    \cT^\lambda_\glo
	    = \sum_{\pm}\sum_{\pm} T^\lambda_{\pm\pm},
	\end{equation}
where the kernel of \(T^\lambda_{++}\) is supported in \(\R_{>0}^2\). The other three choices of signs correspond to the remaining quadrants determined by the signs of \(x\) and \(y\). We henceforth restrict attention to the positive quadrant and write \(T^\lambda_{++}\) simply as \(T^\lambda\); the other three cases are treated by parallel arguments. Accordingly, in \((x,y)\)-coordinates, the relation \(D \simeq D_B\) and all subsequent regions are understood within the positive quadrant \(\R_{>0}^2\).

On the positive quadrant, we now exploit the quasi-homogeneity of \(P\): the Taylor support \(\Lambda(P)\) is contained in an affine line. Since \(\Lambda(P_{xy})\) inherits this collinearity, there exist \(w_1, w_2, d \in \Z\) with \(w_1\) and \(w_2\) not both zero, such that
	\begin{equation} \label{eq:weight_degree}
		w_1m_1+w_2m_2=d
		\quad \text{for every }(m_1,m_2)\in\Lambda(P_{xy}).
	\end{equation}
By passing to the adjoint if necessary, we may assume that \(w_1 \neq 0 \) and set \(\gamma=w_2/w_1 \in \Q\).\footnote{If \(w_1=0\) in \eqref{eq:weight_degree}, then \(P_{xy}(x,y)=y^NQ(x)\), which factors in \(x\) rather than \(y\). In this case, we simply pass to the adjoint operator \((T^\lambda)^*\) to interchange the roles of \(x\) and \(y\).
} Then we have the factorization
	\begin{equation} \label{eq:quasi_form}
		P_{xy}(x,y) = u x^M y^N \prod_{\ell=1}^{L}(y-c_\ell x^\gamma)^{k_\ell} 
		\quad \text{for } x>0.
	\end{equation}
Here \(u \in \R\setminus\set{0}\), \(M,N \in \Z_{\ge 0}\), \(k_\ell \in \N\), and \(c_1, \dotsc, c_L \in \C\setminus\set{0}\) are distinct (see \Cref{lemma:form_quasi} for a self-contained proof). We order the factors so that, for some \(0 \le K \le L\),
	\[
		0<c_1<\dotsb<c_K, \quad \text{and} \quad
		c_\ell \notin \R_{>0} \quad \text{if } K<\ell\le L.
	\]

The factorization \eqref{eq:quasi_form} identifies the real zero curves \(y=c_\ell x^\gamma\) in the positive quadrant. We distinguish the region away from these curves from the regions near them. Since the Hessian and box-area scales differ among these regions, we encode each region by a Newton pair.

\Def{Newton pairs}{\label{def:newton_pairs}
Assume that \(P_{xy}\) is factored as in \eqref{eq:quasi_form}. For \(1 \le \ell \le K\), set \(Y=y-c_\ell x^\gamma\). We define the exponent sets \(\Lambda_0\) and \(\Lambda_\ell\) by
	\begin{equation} \label{eq:Lambda_sets}
		\Lambda_0 \defeq \Lambda(P_{xy}), \quad \text{and} \quad
		\Lambda_\ell \defeq \Lambda\paren[\big]{P_{xy}(x,Y+c_\ell x^\gamma)} 
		\quad \text{for } 1 \le \ell \le K.
	\end{equation}
Here \(\Lambda_\ell\) is the exponent set of \(P_{xy}(x,Y+c_\ell x^\gamma)\), viewed as a finite fractional power series in the \emph{adapted coordinates} \((x,Y)\). We define
	\begin{equation} \label{eq:NP_away}
		\P_0 \defeq \conv\paren{\Lambda_0},
		\qquad \bv_0 \defeq (1,1),
	\end{equation}	
and for \(1 \le \ell \le K\),
	\begin{equation} \label{eq:NP_near}
		\P_\ell \defeq \conv\paren[\big]{\Lambda_\ell+\cone\set{(-\gamma,1)}},
		\qquad
		\bv_\ell \defeq
		\begin{cases}
			(1-\gamma,2), & \text{if } \gamma \neq 0, \\
			(1,1), & \text{if } \gamma =0.
		\end{cases}
	\end{equation}
We call each ordered pair \((\P_\ell, \bv_\ell)\) a \emph{Newton pair} associated with \(P\), and denote the finite family of these pairs by
	\begin{equation} \label{eq:NF}
		\fN(P) \defeq \set[\big]{ (\P_\ell, \bv_\ell) : 0 \le \ell \le K}.
	\end{equation}
When \(P_{xy}\) is a monomial, we have \(L=0\) and the choice of \(\gamma\) is immaterial; in this case, \(\fN(P)\) consists solely of \((\P_0,\bv_0)\).
}
\Rmk{}{
The two components of a Newton pair correspond to the two estimates on a box:
	\begin{center}
		an oscillatory estimate of order \((\lambda|P_{xy}|)^{-1/2}\) and a size estimate of order \((\Delta x \Delta y)^{1/2}\).
	\end{center}
The polygon \(\P_\ell\) determines the Hessian scale \(|P_{xy}|\), while the vector \(\bv_\ell\) records the corresponding box-area scale \(\Delta x \Delta y\).\footnote{
		When \(\gamma\neq0\), one has \(\Delta y \approx x^{\gamma-1}\Delta x\) with \(\Delta y \approx |Y|\) near the curve \(y=c_\ell x^\gamma\). This gives \(\Delta x \Delta y \approx x^{1-\gamma}Y^2\), which accounts for \(\bv_\ell=(1-\gamma,2)\). 
		When \(\gamma=0\), the zero curve is horizontal, and no additional subdivision in \(x\) is needed.
}
The pair \((\P_\ell, \bv_\ell)\) encodes the balance between these two estimates. The precise geometric derivation of these pairs is given in \Cref{sec:preliminaries}.

We emphasize that the vector \(\bv_\ell\) is not an amplitude-weight vector. For comparison, in the weighted scalar oscillatory integral 
	\[
		\iint \ee^{\ii \lambda P(x,y)}x^ay^b \psi(x,y) \dx \dy \quad \text{where } a,b \in \Z_{\ge0},
	\]
the relevant Newton distance is measured along the shifted vector \((1+a,1+b)\), which reflects the prescribed monomial weight together with the integration measure (see, e.g., \cite{KN16}). By contrast, \(\bv_\ell\) arises from the geometry of the operator boxes and records their area scale.
}

We first distinguish balanced and unbalanced Newton pairs according to whether the polygon \(\P_\ell\) meets the ray \(\cone(\bv_\ell)\). For a balanced pair, the endpoints of this intersection define the forward and backward Newton distances.

\Def{Balancedness}{\label{def:balance}
A pair \((\P_\ell,\bv_\ell)\in \fN(P)\) is \emph{balanced} if
	\[
		\P_\ell \cap \cone(\bv_\ell)\neq \emptyset,
	\]
and \emph{unbalanced} otherwise.\footnote{Our balancedness condition is weaker than the balancedness condition in \cite[Definition~2.4]{Kim23}.}
The family \(\fN(P)\) is \emph{balanced} if every pair in \(\fN(P)\) is balanced, and \emph{unbalanced} if at least one pair is unbalanced.
}

\Def{Newton distances}{\label{def:distance}
For a balanced pair \((\P_\ell, \bv_\ell)\), we define the \emph{forward and backward Newton distances} \(\df(\P_\ell, \bv_\ell)\) and \(\db(\P_\ell, \bv_\ell)\) to be the numbers satisfying
	\begin{equation} \label{eq:newton_distance}
		\P_\ell\cap \cone(\bv_\ell)
		= \bracket[\big]{\df(\P_\ell,\bv_\ell),\db(\P_\ell,\bv_\ell)} \bv_\ell.
	\end{equation}
Here, we adopt the convention \(\db(\P_\ell,\bv_\ell)=+\infty\) when the intersection is unbounded.

If the family \(\fN(P)\) of Newton pairs is balanced, we define the \emph{global Newton distances} as
	\begin{equation} \label{eq:global_newton_distances}
		\df(P) \defeq \max_{(\P_\ell, \bv_\ell) \in \fN(P)} \df(\P_\ell,\bv_\ell), \qquad 
		\db(P) \defeq \min_{(\P_\ell, \bv_\ell) \in \fN(P)} \db(\P_\ell,\bv_\ell).
	\end{equation}
}

With these geometric notions in hand, we now state our main theorem.

\begin{main}\label{thm:main}
Let \(P\in \R[x,y]\) be a quasi-homogeneous polynomial and assume \(P_{xy} \not \equiv 0\). Let \(T^\lambda\) denote the positive-quadrant piece in \eqref{eq:T_quadrant_decom}.  Then the following dichotomy holds:
\begin{enumerate}[label=\rm (\Alph*)]
	\item
		If the family \(\fN(P)\) of Newton pairs is balanced, then \(T^\lambda\) extends to a bounded operator on \(L^2(\R)\) for every \(\lambda>0\). Moreover, there exists a constant \(C_P>0\), depending only on \(P\), such that
			\[
				\norm{T^\lambda}_{\op} \le 
				\begin{cases}
					C_P \lambda^{-\frac{1}{2(\df(P)+1)}}, & \text{if } \lambda  \ge 1, \\
					C_P \lambda^{-\frac{1}{2(\db(P)+1)}}, & \text{if } 0< \lambda < 1.
				\end{cases}
			\]
		Here, \(\df(P)\) and \(\db(P)\) are defined as in \eqref{eq:global_newton_distances}.
		These bounds are sharp: the corresponding reverse inequalities hold with a positive constant depending only on \(P\).

	\item
		If the family \(\fN(P)\) of Newton pairs is unbalanced, then for every \(\lambda>0\),
			\[
				\sup_{0 \neq f \in \cS(\R)} \frac{\norm{T^\lambda f}_2}{\norm{f}_2} = +\infty.
			\]
\end{enumerate}
\end{main}

Since each quadrant piece is a spatial restriction of \(\cT^\lambda_\glo\), one has \(\norm{T^\lambda_{\pm\pm}}_{\op} \le \norm{\cT^\lambda_\glo}_{\op}\). Applying \Cref{thm:main} to each quadrant piece in \eqref{eq:T_quadrant_decom}, we completely resolve the global problem for quasi-homogeneous polynomial phases.

\Cor{}{\label{cor:full_operator}
For the operator \(\cT^\lambda_\glo\) defined in \eqref{eq:T_glo}, the following dichotomy holds:
\begin{enumerate}[label=\rm (\Alph*)]
	\item
		If every Newton pair arising from the four quadrant pieces is balanced, then \(\cT^\lambda_\glo\) extends to a bounded operator on \(L^2(\R)\) for every \(\lambda>0\). Moreover, the sharp estimates
			\[
				\norm{\cT^\lambda_\glo}_{\op} \approx 
				\begin{cases}
					\lambda^{-\frac{1}{2(d_1+1)}}, & \text{if } \lambda  \ge 1, \\
					\lambda^{-\frac{1}{2(d_2+1)}}, & \text{if } 0< \lambda < 1,
				\end{cases}
			\]
		hold, where \(d_1\) and \(d_2\) are the largest forward and the smallest backward Newton distances among these pairs, respectively.

	\item
		If any quadrant contains an unbalanced Newton pair, then \(\cT^\lambda_\glo\) has no bounded extension to \(L^2(\R)\) for any \(\lambda>0\).
\end{enumerate}
}

\Rmk{}{\label{rmk:explicit_distances}
Remarkably, thanks to quasi-homogeneity, the distances \(d_1\) and \(d_2\) of \Cref{cor:full_operator}\,(A) can be computed purely algebraically from the factorization of \(P_{xy}\). If we write \(P_{xy}\) as in \eqref{eq:quasi_form},
	\[
		P_{xy}(x,y) = u x^M y^N \prod_{\ell=1}^{L}(y-c_\ell x^\gamma)^{k_\ell} 
		\quad \text{for } x>0,
	\]
then \(d_1\) and \(d_2\) are given explicitly:
	\[
		\begin{alignedat}{2}
			&\text{if }\gamma \neq -1, \qquad
				&&d_1=d_2
					=\frac{M+\gamma\paren[\big]{N+\sum_{\ell=1}^L k_\ell}}{1+\gamma},\\[0.5ex]
			&\text{if } \gamma = -1, \qquad
				&&d_1=\max\set[\Big]{N, \, \max_{c_\ell \in \R} \frac{k_\ell}{2}},
					\ d_2=M=N+\sum_{\ell=1}^L k_\ell.
		\end{alignedat}
	\]
These formulas follow directly from the collinearity of the Newton polygons forced by quasi-homogeneity. Here the equality \(M=N+\sum_{\ell=1}^L k_\ell\) in the second line follows from balancedness. When no \(c_\ell\) is real, the inner maximum is omitted. Moreover, sign changes \((x,y) \mapsto (\pm x, \pm y)\) preserve \(M,N,\gamma\) and \(\sum_{\ell=1}^L k_\ell\), and, when \(\gamma=-1\), send each real root \(c_\ell\) to \(\pm c_\ell\) with the same multiplicity \(k_\ell\), so that \(\max_{c_\ell \in \R} k_\ell/2\) is unchanged. Consequently, the algebraic expressions above do not depend on the choice of quadrant.
}

\begin{example*}
For the phase \(P(x,y)=-x^2y\) in \Cref{ex:1@intro}, one has \(P_{xy}(x,y)=-2x\). Hence every quadrant yields the unique Newton pair \((\P_0, \bv_0)\) where
	\[
		\P_0 = \set{(1,0)}, \qquad \bv_0=(1,1).
	\]
Since \(\P_0 \cap \cone(\bv_0) = \emptyset\), this pair is unbalanced, and \Cref{cor:full_operator}\,(B) recovers the global \(L^2\)-unboundedness established directly in that example.

Moreover, for the phase \(P(x,y)=x^2y^2-4xy\) from \Cref{ex:2@intro,ex:3@intro}, one has
	\[
		P_{xy}(x,y)=4x(y-x^{-1}), \qquad 
		P_{xy}(x,Y+x^{-1}) = 4xY \quad \text{for } x,y>0.
	\]
Thus we have
	\[
		\begin{alignedat}{2}
			&\P_0 = \conv\set{(0,0), (1,1)}, \qquad 
				&& \P_1 = (1,1)+ \cone\set{(1,1)}, \\
			&\bv_0 = (1,1), 
				&& \bv_1 = (2,2).
		\end{alignedat}
	\]
Consequently, \(\df(\P_0, \bv_0)=0\), \(\db(\P_0,\bv_0)=1\) and \(\df(\P_1, \bv_1) = 1/2\), \(\db(\P_1, \bv_1)=+\infty\). The same Newton pairs arise from the third quadrant, while the second and fourth quadrants contain no zeros of \(P_{xy}\) and contribute only \((\P_0, \bv_0)\). The global Newton distances are therefore \(d_1=\max\set{0,1/2}=1/2\) and \(d_2=\min\set{1,+\infty}=1\), which exactly recover the sharp global rates \(\lambda^{-1/3}\) and \(\lambda^{-1/4}\).
\end{example*}

\subsection{Regimes where a scaling argument applies or fails} \label{subsec:scaling}

By applying a scaling argument to the operator \(\cT^\lambda_\glo\) in \eqref{eq:T_glo}, we now clarify why the global problem is not a formal consequence of the local theory. By \eqref{eq:weight_degree}, one has \(P_{xy}(R^{w_1}x, R^{w_2}y) = R^{d} P_{xy}(x,y)\) for any \(R>0\). Since \(P(x,0)=P(0,y)=0\), this relation integrates to
	\begin{equation} \label{eq:scaling_P}
		P(R^{w_1}x, R^{w_2}y) = R^{d+w_1+w_2}P(x,y).
	\end{equation}
Let \(\psi \in \cC_0^\infty(\R)\) be supported sufficiently close to the origin and satisfy \(\psi \equiv 1\) in a neighborhood of the origin. For fixed \(\lambda >0\), we define
	\begin{equation} \label{eq:T_trunc}
		\cT^\lambda_R f(x) \defeq
		\int \ee^{\ii\lambda P(x,y)}
			\psi\paren[\Big]{\frac{x}{R^{w_1}}}\psi\paren[\Big]{\frac{y}{R^{w_2}}}f(y)\dy.
	\end{equation}
Then the scaling relation \eqref{eq:scaling_P} and the change of variables \((x,y) \mapsto (R^{w_1}x, R^{w_2}y)\) yield
	\[
		\norm{\cT^\lambda_R}_{\op}
		= R^{\frac{w_1+w_2}{2}}\norm{\cT_1^{\lambda R^{d+w_1+w_2}}}_{\op},
	\]
where \(\cT_1^{\widetilde \lambda}\) denotes the operator in \eqref{eq:T_trunc} with \(R=1\) and parameter \(\widetilde \lambda\). Applying \Cref{thm:PS@intro} to the operator \(\cT^{\lambda R^{d+w_1+w_2}}_1\), we obtain
	\begin{align}
		\norm{\cT^\lambda_R}_{\op} 
		&\approx R^{\frac{w_1+w_2}{2}} (\lambda R^{d+w_1+w_2})^{-\frac{1}{2(\delta_\loc+1)}} 
			\notag \\
		&= \lambda^{-\frac{1}{2(\delta_\loc+1)}}R^{\frac{(w_1+w_2)\delta_\loc-d}{2(\delta_\loc+1)}} \quad \text{for } \lambda R^{d+w_1+w_2} \ge 1. 
			\label{eq:dilation_block_norm}
	\end{align}

\paragraph{The case \(w_1w_2>0\).}
By replacing \((w_1, w_2, d)\) with \((-w_1, -w_2, -d)\) in \eqref{eq:weight_degree} if necessary, we may further assume that \(w_1, w_2>0\). The supports of the cutoffs in \eqref{eq:T_trunc} then exhaust \(\R^2\) as \(R \toinf\).
Recall from \labelcref{eq:local_NP,eq:local_newton_distance} that the local Newton polygon and distance are defined by
	\[
		\bN(P_{xy}) = \conv\paren[\big]{\Lambda(P_{xy})+\cone\set{\be_1,\be_2}},
		\qquad
		\bN(P_{xy})\cap \cone(\1) = [\delta_{\loc},\infty)\1.
	\]
By \eqref{eq:weight_degree}, every point of \(\conv(\Lambda(P_{xy}))\) lies on the line \(w_1m_1+w_2m_2=d\), which intersects the diagonal ray \(\cone(\1)\) at the point \(\frac{d}{w_1+w_2}\1\).
Since \(w_1, w_2>0\), it follows that
	\begin{enumerate}[label=(\roman*)]
		\item 
			if \(\conv(\Lambda(P_{xy})) \cap \cone(\1) \neq \emptyset\), then \((w_1+w_2)\delta_\loc =d\);
		\item
			if \(\conv(\Lambda(P_{xy}))\cap \cone(\1) = \emptyset\), then \((w_1+w_2)\delta_\loc > d\).
	\end{enumerate}
Therefore, in case (i), substituting \((w_1+w_2)\delta_\loc = d\) into \eqref{eq:dilation_block_norm} and letting \(R\toinf\) yields
	\[
		\norm{\cT^\lambda_\glo}_{\op} 
		\approx \lambda^{-\frac{1}{2(\delta_\loc+1)}}.
	\]
This local Newton distance \(\delta_\loc\) exactly coincides with both global Newton distances in \Cref{cor:full_operator}: direct computation yields \(\delta_\loc  = d/(w_1+w_2)=d_1=d_2\).
In case (ii), we have \(\norm{\cT^\lambda_R}_{\op} \toinf\) as \(R\toinf\) by \eqref{eq:dilation_block_norm}. Hence \(\cT^\lambda_\glo\) has no bounded extension to \(L^2(\R)\). Thus the local theory combined with the scaling argument fully resolves the regime \(w_1w_2>0\).

Phong and Stein \cite[Theorem~1 and Remark~2.1]{PS94} already obtained the complete global \(L^2\)-classification for homogeneous models: their arguments yield the sharp decay in the bounded case and \(L^2\)-unboundedness otherwise. The same arguments readily extend to quasi-homogeneous phases with \(w_1w_2>0\).

\paragraph{The case \(w_1w_2 \le 0\).}
In this regime, the supports of the cutoffs in \eqref{eq:T_trunc} no longer exhaust \(\R^2\). Consequently, the scaling argument cannot reduce the global problem to the local model. This failure is already visible in \Cref{ex:2@intro}: for \(P(x,y)=x^2y^2-4xy\), the sharp global decay rate is strictly slower than its local counterpart. Phong and Stein left the complementary cases \(w_1 w_2 \le 0\) unresolved. We completely resolve these remaining cases.


\subsection{An application to a global Radon transform}

As a direct application of \Cref{cor:full_operator}, we determine the sharp anisotropic Sobolev regularity for a Radon transform. In the compactly supported setting, Phong and Stein \cite{PS94} considered local Radon transforms of the form
	\[
		\cR_\loc g(t,x) = \int g\paren[\big]{t+P(x,y),y}\psi(t,x,y)\dy,
	\]
where \(\psi \in \cC_0^\infty(\R^3)\) is a smooth cutoff function. Motivated by this local model, we define the corresponding global Radon transform for \(g=g(s,y) \in \cS(\R^2)\) by
	\begin{equation} \label{eq:global_Radon}
		\cR g(t,x) \defeq \int g\paren[\big]{t+P(x,y), y} \dy.
	\end{equation}

The partial Fourier transform in \(t\) yields a fiber decomposition of \(\cR\), with fiber \(\cT^\lambda_\glo\) at each frequency \(\lambda \neq 0\).  Consequently, \Cref{cor:full_operator} determines the exact range of \(\tau \in \R\) for which \(\cR\) extends to a bounded operator from \(L^2(\R^2)\) to the anisotropic homogeneous Sobolev space \(\dot H^{(\tau,0)}(\R^2)\), equipped with the seminorm
	\[
		\norm{F}_{\dot H^{(\tau,0)}(\R^2)} 
		\defeq \norm[\big]{|D_t|^\tau F}_{L^2(\R^2)}.
	\]
	
In the statement below, \(d_1\) and \(d_2\) are as in \Cref{cor:full_operator}. Since \(P\) is real, the conclusions of \Cref{cor:full_operator} apply equally to negative \(\lambda\), and hence control all Fourier frequencies in the \(t\)-variable.

\Cor{Anisotropic Sobolev regularity}
{\label{cor:global_Radon}
Let \(P \in \R[x,y]\) be a quasi-homogeneous polynomial with \(P_{xy} \not \equiv 0\), and let \(\cR\) be defined by \eqref{eq:global_Radon}. Then the following hold.
\begin{enumerate}[label=\rm (\Alph*)]
	\item
		If every Newton pair arising from the four quadrant pieces is balanced, then for \(\tau \in \R\),
			\[
				\norm[\big]{|D_t|^\tau \cR g}_{L^2(\R^2)} \lesssim \norm{g}_{L^2(\R^2)} 
				\quad \text{for all } g \in \cS(\R^2)
			\]
		holds if and only if
			\[
				\frac{1}{2(d_2+1)} \le \tau \le \frac{1}{2(d_1+1)}.
			\]
		Equivalently, \(\cR\) extends to a bounded operator from \(L^2(\R^2)\) to \(\dot H^{(\tau,0)}(\R^2)\) precisely for this range of \(\tau\).
		
	\item
		If any quadrant contains an unbalanced Newton pair, then \(\cR\) has no bounded extension from \(L^2(\R^2)\) to \(\dot H^{(\tau,0)}(\R^2)\) for any \(\tau \in \R\).
\end{enumerate}
}

The proof of \Cref{cor:global_Radon} is standard and follows from Plancherel's theorem; we defer the details to \Cref{appendix:radon}.

The range of \(\tau\) in \Cref{cor:global_Radon}\,(A) is always nonempty. Indeed, the explicit formulas in \Cref{rmk:explicit_distances} show that when \(\gamma \neq -1\), one has \(d_1=d_2\) and the \(\tau\)-interval collapses to a single point. When \(\gamma =-1\), the same formulas give
	\[
		d_1=\max\set[\Big]{N, \, \max_{c_\ell \in \R} \frac{k_\ell}{2}} \le N+\sum_{\ell=1}^L k_\ell = M = d_2.
	\]
In the latter case, the interval is nontrivial precisely when \(d_1 < d_2\), that is, when the sharp decay rates in the high- and low-frequency regimes are different.

In the local setting, Phong and Stein \cite[Theorems~3 and 4]{PS94} established both \(L^p\)-\(L^q\) improving estimates and isotropic Sobolev regularity. In the global setting, however, isotropic Sobolev regularity fails in general. As shown in \Cref{appendix:radon}, for the nondegenerate phase \(P(x,y)=xy\), the transform \(\cR\) admits no bounded estimate with nonzero Sobolev order in the \(x\)-variable. Thus the restriction to \(t\)-regularity in \Cref{cor:global_Radon} is intrinsic to the cutoff-free problem.


\subsection{Further related work}

\paragraph{Local oscillatory integral operators.}
Phong and Stein first treated polynomial phases and homogeneous models \cite{PS92,PS94} before developing their local Newton-polygon theory \cite{PS97}. Rychkov \cite{Ryc01} extended the result of Phong and Stein to \(\cC^\infty\) phases, with a possible loss of a power of \(\log \lambda\). His proof combines Puiseux analysis of the real branches of \(\set{S_{xy}=0}\), dyadic summation away from them, and anisotropic Whitney decomposition near them. The latter decomposition is especially close to the admissible-box construction used in the present paper. Greenblatt \cite{Gre04,Gre05} subsequently developed a direct resolution procedure tailored to oscillatory integral operators and removed the logarithmic loss, obtaining sharp \(L^2\)-estimates for \(\cC^\infty\) phases.

\paragraph{Scalar oscillatory integrals.}
On the scalar Newton-polyhedral side, Varchenko \cite{Var76} studied oscillatory decay under a suitable non-degeneracy condition, while Vassiliev \cite{Vas77} developed related Newton-diagram methods for exponential-integral asymptotics. Phong, Stein, and Sturm \cite{PSS99,PSS01} applied Newton-polyhedral methods to growth, stability, and multilinear level-set problems. Ikromov and Müller \cite{IM11,IM16} developed adapted coordinate systems and principal root jets in two-dimensional Newton geometry; their notion of adaptedness, which maximizes the Newton distance over coordinate changes, is unrelated to the adapted coordinates \((x,Y)\) of the present paper, which are simply shears along a fixed real zero curve. Greenblatt \cite{Gre10} related scalar oscillatory decay and sublevel-set growth to Newton polyhedra, while Collins, Greenleaf, and Pramanik \cite{CGP13} developed a higher-dimensional resolution procedure.

\paragraph{Degenerate canonical relations.}
For canonical relations whose two projections are folds, Pan and Sogge \cite{PanSog90} obtained the characteristic \(1/6\)-derivative loss, while Greenleaf and Seeger \cite{GS94} established sharp estimates when one projection is a fold, with no hypothesis on the other projection. Seeger \cite{See93} studied finite-type degeneracies in the plane. Greenleaf and Seeger later treated cusp singularities and further low-order corank-one degeneracies, including swallowtails \cite{GS98,GS02b}, and surveyed degenerate canonical relations \cite{GS02a}. Cuccagna \cite{Cuc97} obtained \(L^2\)-estimates for averaging operators with two-sided \(k\)-fold singularities, while Christ, Nagel, Stein, and Wainger \cite{CNSW99} analyzed singular and maximal Radon transforms under finite-type hypotheses.

\paragraph{\(L^p\)-estimates and higher dimensions.}
Beyond the local \(L^2\)-theory, Yang \cite{Yan04} obtained sharp one-dimensional \(L^p\)-estimates for oscillatory integral operators with polynomial phases, and Xiao \cite{Xia17} established endpoint estimates for one-dimensional oscillatory integral operators. Shi, Xu, and Yan \cite{SXY19} developed damping estimates for the local analytic \(L^p\)-theory. Tang \cite{Tan06} studied decay rates for homogeneous phases in \(1+2\) dimensions under a non-degeneracy assumption. In the compactly supported setting, Greenleaf, Pramanik, and Tang \cite{GPT07} treated homogeneous polynomial phases in several variables at the \(L^2\) level under genericity assumptions, while Xu and Yan \cite{XY19} obtained higher-dimensional \(L^p\)-estimates formulated in terms of Newton distance. In a recent preprint, Xu \cite{Xu26} obtained the sharp local \(L^2\)-decay for \((2+1)\)-dimensional operators with phase \(P(x,y)t^k\), governed by the competition between Varchenko's Newton height of \(P\) and the degeneracy of \(t^k\), and also derived a sharp \(L^2\to L^{2k+2}\) estimate.

\paragraph{Newton-polyhedral methods on unbounded domains.}
The second author previously used faces, cones, and cone-type decompositions of Newton polyhedra to characterize the boundedness of multiple Hilbert transforms associated with polynomial mappings \cite{Kim15}, and subsequently developed a domain-adapted Newton-polyhedral framework, including balancing conditions and forward and backward distances, for scalar sublevel sets and oscillatory integrals over unbounded domains \cite{Kim23}. The present paper adapts this convex-geometric perspective to the cutoff-free \(L^2\)-operator problem. The operator-specific ingredients are the Newton pair \((\P_\ell, \bv_\ell)\) and the admissible-box construction, which makes the Hessian bounds obtained in the adapted coordinates available for operator estimates in the original coordinates.


\subsection*{Idea of the proof and organization}

\paragraph{Idea of the proof.}
The principal geometric difficulty in estimating the global operator is that the Hessian geometry and the operator estimates naturally belong to different coordinate systems. More precisely, near a real zero curve \(y=c_\ell x^\gamma\), the \emph{adapted} coordinates \((x,Y)\), where \(Y=y-c_\ell x^\gamma\), put the Hessian into a normal-crossing form. The resulting expansion of \(P_{xy}(x,Y+c_\ell x^\gamma)\) in \((x,Y)\) determines the adapted Newton polygon \(\P_\ell\). However, all operator estimates are tied to the original coordinate structure. In particular, the operator van der Corput estimate is applied in the original coordinates, and the size estimate and the almost-orthogonality argument depend on the projections of the kernel support onto the original \(x\)- and \(y\)-axes. We must therefore organize the singular geometry in the \emph{adapted} coordinates while carrying out the operator estimates in the \emph{original} coordinates. Greenblatt \cite{Gre04} explicitly emphasized the necessity of this separation. Phong, Stein, and Sturm \cite{PSS01} also discussed the coordinate dependence underlying this separation from a broader perspective.

We overcome this difficulty by connecting the two coordinate systems through admissible boxes. These boxes are axis-parallel in the \emph{original} \((x,y)\)-coordinates and retain the Hessian bounds obtained from the \emph{adapted} expansion \(P_{xy}(x,Y+c_\ell x^\gamma)\). They capture the same scales and overlap properties as the anisotropic Whitney rectangles used by Rychkov \cite{Ryc01}; we encode the area scale of the admissible boxes by the exponent vector \(\bv_\ell\). Thus, \(\P_\ell\) determines the Hessian scale entering the oscillatory estimate, whereas \(\bv_\ell\) records the area scale entering the size estimate. 

We combine these two geometric objects into the Newton pair \((\P_\ell,\bv_\ell)\). The global oscillatory integral operator is \(L^2\)-bounded if and only if every Newton pair arising from the four quadrants is balanced: each \(\P_\ell\) intersects the corresponding ray \(\cone(\bv_\ell)\). In the balanced case, the corresponding Newton distances determine the sharp decay rates. As the preceding examples demonstrate, the local Newton distance alone neither characterizes global \(L^2\)-boundedness nor determines the sharp global decay rate. The Newton pairs therefore provide the geometric framework required to establish the complete global \(L^2\)-theory.

\paragraph{Organization of the paper.}
This paper is organized as follows. The geometric decomposition, admissible boxes, and analytic estimates are developed in \Cref{sec:preliminaries}. We establish the upper bound, the matching lower bound, and the \(L^2\)-unboundedness result in \Cref{sec:upper_bounds,sec:lower_bounds,sec:unboundedness}, respectively. In \Cref{sec:conclusion}, we provide a convex-geometric summary and discuss the extension to general polynomial phases. Finally, we prove the factorization lemma for quasi-homogeneous polynomials in \Cref{appendix:form_quasi} and give the details of the Sobolev-regularity application to the global Radon transform in \Cref{appendix:radon}.

\paragraph{Notation.}
We write \(\R_{\ge0}\) and \(\Z_{\ge0}\) for the sets of nonnegative real numbers and nonnegative integers, respectively, and \(\R_{>0}\) for the set of positive real numbers.

We employ the following auxiliary cutoff functions:
	\begin{itemize}
		\item
			\(\psi \in \cC_0^\infty(\R)\) such that \(\supp \psi \subset [-1,1]\),  \(0 \le \psi \le 1\) and \(\psi \equiv 1\) on \([-1/2, 1/2]\);
		\item 
			\(\chi(t)=\psi(t/2)-\psi(t)\).
	\end{itemize}
We allow slight modifications of \(\psi\) and \(\chi\) from line to line.

For nonnegative quantities \(a\) and \(b\), we write \(a \lesssim b\) if \(a\le Cb\) for some constant \(C>0\). Unless otherwise stated, the implicit constants may depend on the phase \(P\) and the auxiliary cutoff functions, but are independent of the spatial variables \(x,y\), the parameter \(\lambda\), and all dyadic indices. The implicit constants may change from line to line. We write \(a \gtrsim b\) if \(b \lesssim a\), and \(a \approx b\) if both \(a \lesssim b\) and \(b \lesssim a\).


\section{Preliminaries} \label{sec:preliminaries}

We keep the notation introduced in \labelcref{eq:quasi_form,eq:Lambda_sets,eq:NP_away,eq:NP_near}. We first derive the Newton polygons and construct the admissible boxes. We then analyze the geometry of these polygons and record the Hessian bounds and analytic estimates needed in the proofs.


\subsection{Operator decomposition and admissible boxes} \label{subsec:decomposition}

The case in which \(P_{xy}\) has no zeros in \(\R_{>0}^2\) is covered directly by the subsequent analysis away from the zero set. We therefore focus on the case where \(P_{xy}\) has real zero curves in \(\R_{>0}^2\).

By the factorization 
	\[
		P_{xy}(x,y)= ux^My^N\prod_{\ell=1}^L (y-c_\ell x^\gamma)^{k_\ell} \quad \text{for } x>0,
	\]
these real zero curves are given by \(y=c_\ell x^\gamma\) for \(1\le \ell \le K\). Recall that \(0< c_1 < \dotsb< c_K\). To isolate these zero curves, we fix a constant \(\epsilon_0>0\), depending only on \(P\) and the auxiliary cutoffs, such that
	\begin{equation} \label{eq:epsilon_separation}
		0< \epsilon_0 < 2^{-(5|\gamma|+5)} 
			\min \set[\big]{1, \, c_1, \, \min_{1\le m < n \le L} |c_m-c_n|}.
	\end{equation}
When \(L=1\), we omit the inner minimum. When \(\gamma \neq 0\), we additionally require
	\begin{equation} \label{eq:epsilon_box}
		0< \epsilon_0 \le 2^{-(5|\gamma|+5)}c_1|\gamma|.
	\end{equation}
For this choice of \(\epsilon_0\), the regions
	\[
		\set{(x,y) : x>0, \ |y- c_\ell x^\gamma| \le \epsilon_0 x^\gamma}
		\quad \text{for } \ell =1, \dotsc, K,
	\]
are mutually disjoint and contained in \(\R_{>0}^2\). Using the cutoff \(\psi\), we define
	\begin{equation} \label{eq:supp_fts}
		\psi_\ell(x,y) \defeq \psi\paren[\Big]{\frac{y-c_\ell x^\gamma}{\epsilon_0 x^\gamma}}
		\quad \text{and} \quad
		\psi_0(x,y) \defeq 1- \sum_{\ell=1}^K \psi_\ell(x,y).
	\end{equation}
Since \(\psi_0 + \sum_{\ell=1}^K \psi_\ell =1\), we decompose the operator \(T^\lambda\) as \(T^\lambda = \sum_{\ell=0}^K T^\lambda_\ell\), where
	\begin{align}
		T^\lambda_0 f(x) &\defeq \int \ee^{\ii\lambda P(x,y)}\psi_0(x,y) f(y) \dy, 
			\label{eq:T_decom_away} \\ 
		T^\lambda_\ell f(x) &\defeq \int \ee^{\ii\lambda P(x,y)}\psi_\ell(x,y) f(y) \dy 
		\quad \text{for } 1 \le \ell \le K.  
			\label{eq:T_decom_near}
	\end{align}

\paragraph{The piece \(T^\lambda_0\).}
We first consider \(T^\lambda_0\) in \eqref{eq:T_decom_away}. The cutoff \(\psi_0\) in \eqref{eq:supp_fts} removes only small neighborhoods of the real zero curves and places no restriction on the dyadic sizes of \(x\) and \(y\). In the notation of \eqref{eq:domain}, this absence of size restrictions is encoded by \(D_{\set{\bm 0}}=\R^2\), as in \Cref{ex:NP}. Therefore, by \eqref{eq:NP}, we assign to \(T^\lambda_0\) the Newton polygon
	\[
		\bN\paren[\big]{P_{xy}(x,y), D_{\set{\bm 0}}} = \conv(\Lambda_0),
	\]
which is precisely \(\P_0\) in \eqref{eq:NP_away}. 

To identify the vector \(\bv_0\) in \eqref{eq:NP_away}, we decompose \(T^\lambda_0\) into ordinary dyadic blocks in the original coordinates:
	\[
		T^\lambda_0 = \sum_{J=(j_1, j_2) \in \Z^2} T^\lambda_{0,J},
		\quad \text{where }
		T^\lambda_{0,J} f(x) 
			= \int \ee^{\ii\lambda P(x,y)}\psi_0(x,y)
				\chi\paren[\Big]{\frac{x}{2^{-j_1}}}\chi\paren[\Big]{\frac{y}{2^{-j_2}}}f(y)\dy.
	\]
The kernel of \(T^\lambda_{0,J}\) is supported in an axis-parallel box on which \(x \approx 2^{-j_1}\) and \(y \approx 2^{-j_2}\). The area scale of this block is therefore given by
	\[
		2^{-j_1}2^{-j_2} \approx x^1y^1.
	\]
In view of the exponent vector \((1,1)\) of \(x^1y^1\), this justifies the choice of \(\bv_0 = (1,1)\) in \eqref{eq:NP_away}.

\paragraph{The pieces \(T^\lambda_\ell\) \((1\le \ell \le K)\).}
We next fix \(1 \le \ell \le K\) and consider \(T^\lambda_\ell\) in \eqref{eq:T_decom_near}. Setting \(Y=y-c_\ell x^\gamma\), the cutoff \(\psi_\ell\) in \eqref{eq:supp_fts} can be written as \(\psi_\ell(x,y) = \psi(x^{-\gamma}Y/\epsilon_0)\); hence, its support satisfies \(|x|^{-\gamma}|Y|\lesssim 1\). By \Cref{def:cone_domain}, this support condition defines a region in the adapted \((x,Y)\)-coordinates equivalent to \(D_{\set{(-\gamma,1)}}\). Therefore, by \eqref{eq:NP}, we assign to \(T^\lambda_\ell\) the Newton polygon
	\[
		\bN\paren[\big]{P_{xy}(x,Y+c_\ell x^\gamma), D_{\set{(-\gamma,1)}}}
		= \conv\paren[\big]{\Lambda_\ell + \cone\set{(-\gamma,1)}},
	\]
which is precisely \(\P_\ell\) in \eqref{eq:NP_near}.

To identify \(\bv_\ell\) in \eqref{eq:NP_near}, we first localize the \(x\)-scale by
	\[
		T^\lambda_\ell = \sum_{j_1 \in \Z} T^\lambda_{\ell,j_1},
		\quad \text{where }
		T^\lambda_{\ell,j_1} f(x)
			= \int \ee^{\ii\lambda P(x,y)}\psi_\ell(x,y)
				\chi\paren[\Big]{\frac{x}{2^{-j_1}}} f(y) \dy.
	\]
Fix \(j_1\), and set \(I_{j_1}=[2^{-j_1-1}, 2^{-j_1+1}]\).
We work on one side of the curve \(y=c_\ell x^\gamma\); the other side is treated identically.

\subparagraph{Non-horizontal curves.}
We first treat the case \(\gamma\neq0\). In the adapted coordinates \((x,Y)\) where \(Y=y-c_\ell x^\gamma\), the region under consideration is given by \(x \in I_{j_1}\) and \(Y \approx h\). In the original coordinates \((x,y)\), however, it is a curved strip along \(y=c_\ell x^\gamma\).

To carry out the operator estimates in the original coordinates, we subdivide \(I_{j_1}\) until the variation of \(x \mapsto c_\ell x^\gamma\) is \(\approx h\), and cover each resulting curved piece by an axis-parallel box. The following proposition records this construction and the resulting overlap and area estimates.

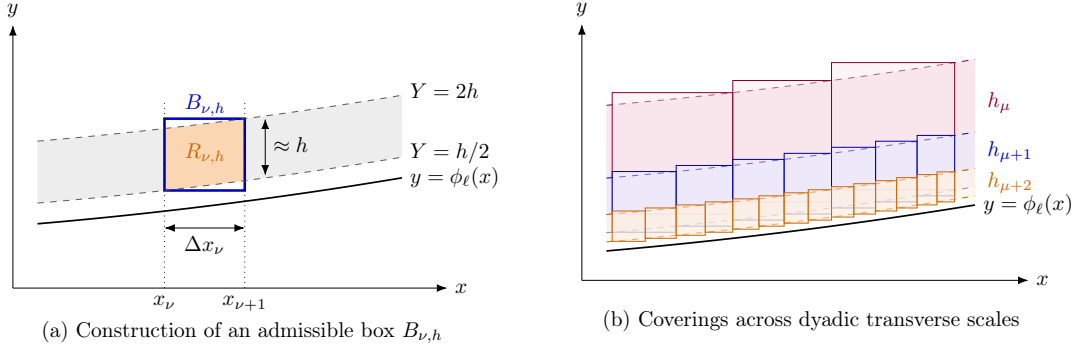
\begin{figure}[htb]
\centering
\resizebox{0.9\textwidth}{!}{%
\begin{tikzpicture}[scale=0.75, >=Latex, font=\small]

\begin{scope}[xshift=0cm, xscale=1.8, yscale=1.65]
	\pgfmathsetmacro{\h}{0.55}
	\pgfmathsetmacro{\xL}{0.65}
	\pgfmathsetmacro{\xR}{5.10}
	\pgfmathsetmacro{\xA}{2.20}
	\pgfmathsetmacro{\xB}{3.18}
	\pgfmathsetmacro{\phiA}{1+0.08*\xA+0.010*\xA*\xA}
	\pgfmathsetmacro{\phiB}{1+0.08*\xB+0.010*\xB*\xB}

	\draw[->] (0.35,0.20) -- (5.65,0.20) node[right] {$x$};
	\draw[->] (0.35,0.20) -- (0.35,3.70) node[above] {$y$};

	\path[fill=gray!14]
		plot[domain=\xL:\xR, samples=100, variable=\x]
			(\x,{1+0.08*\x+0.010*\x*\x+2*\h})
		--
		plot[domain=\xR:\xL, samples=100, variable=\x]
			(\x,{1+0.08*\x+0.010*\x*\x+\h/2})
		-- cycle;

	\path[fill=orange!35, opacity=0.82]
		plot[domain=\xA:\xB, samples=60, variable=\x]
			(\x,{1+0.08*\x+0.010*\x*\x+2*\h})
		--
		plot[domain=\xB:\xA, samples=60, variable=\x]
			(\x,{1+0.08*\x+0.010*\x*\x+\h/2})
		-- cycle;

	\draw[blue!75!black, very thick]
		(\xA,{\phiA+\h/2}) rectangle
		(\xB,{\phiB+2*\h});

	\draw[thick]
		plot[domain=\xL:\xR, samples=100, variable=\x]
			(\x,{1+0.08*\x+0.010*\x*\x});

	\draw[dashed, gray!70!black]
		plot[domain=\xL:\xR, samples=100, variable=\x]
			(\x,{1+0.08*\x+0.010*\x*\x+\h/2});

	\draw[dashed, gray!70!black]
		plot[domain=\xL:\xR, samples=100, variable=\x]
			(\x,{1+0.08*\x+0.010*\x*\x+2*\h});

	\node[right] at
		(5.10,{1+0.08*5.10+0.010*5.10*5.10})
		{$y=\phi_\ell(x)$};

	\node[right, yshift=2pt] at
		(5.10,{1+0.08*5.10+0.010*5.10*5.10+\h/2})
		{$Y=h/2$};

	\node[right, yshift=2pt] at
		(5.10,{1+0.08*5.10+0.010*5.10*5.10+2*\h})
		{$Y=2h$};

	\draw[dotted]
		(\xA,0.20) -- (\xA,{\phiB+2*\h+0.22});
	\draw[dotted]
		(\xB,0.20) -- (\xB,{\phiB+2*\h+0.22});

	\node[below] at (\xA,0.20) {$x_\nu$};
	\node[below] at (\xB,0.20) {$x_{\nu+1}$};

	\draw[<->] (\xA,1) -- (\xB,1)
		node[midway, below] {$\Delta x_\nu$};

	\draw[<->]
		({\xB+0.25},{\phiA+\h/2+0.2}) --
		({\xB+0.25},{\phiB+2*\h})
		node[midway, right, yshift=1pt] {$\approx h$};

	\node[orange!80!black]
		at (2.70,{\phiA+1.25*\h+0.08})
		{$R_{\nu,h}$};

	\node[blue!75!black]
		at (2.70,{\phiB+2*\h+0.18})
		{$B_{\nu,h}$};

	\node at (3.15,-0.38)
		{\small (a) Construction of an admissible box \(B_{\nu,h}\)};
\end{scope}

\begin{scope}[xshift=12.5cm,xscale=1.8,yscale=2.5]
	\pgfmathsetmacro{\xL}{0.65}
	\pgfmathsetmacro{\xR}{5.15}
	\pgfmathsetmacro{\hA}{0.64}
	\pgfmathsetmacro{\hB}{\hA/2}
	\pgfmathsetmacro{\hC}{\hA/4}

	\draw[->] (0.35,0.20) -- (5.80,0.20) node[right] {$x$};
	\draw[->] (0.35,0.20) -- (0.35,2.45) node[above] {$y$};

	\path[fill=purple!10]
		plot[domain=\xL:\xR,samples=100,variable=\x]
			(\x,{0.42+0.055*\x+0.006*\x*\x+2*\hA})
		--
		plot[domain=\xR:\xL,samples=100,variable=\x]
			(\x,{0.42+0.055*\x+0.006*\x*\x+\hA/2})
		-- cycle;
	\draw[dashed,purple!70!black]
		plot[domain=\xL:\xR,samples=100,variable=\x]
			(\x,{0.42+0.055*\x+0.006*\x*\x+\hA/2});
	\draw[dashed,purple!70!black]
		plot[domain=\xL:\xR,samples=100,variable=\x]
			(\x,{0.42+0.055*\x+0.006*\x*\x+2*\hA});
	\foreach \xa/\xb in {0.720/2.192,2.192/3.397,3.397/4.900}{
		\pgfmathsetmacro{\pa}{0.42+0.055*\xa+0.006*\xa*\xa}
		\pgfmathsetmacro{\pb}{0.42+0.055*\xb+0.006*\xb*\xb}
		\draw[purple!70!black,line width=0.45pt,fill=purple!10,fill opacity=0.42]
			(\xa,{\pa+\hA/2}) rectangle (\xb,{\pb+2*\hA});
	}

	\path[fill=blue!8,fill opacity=0.55]
		plot[domain=\xL:\xR,samples=100,variable=\x]
			(\x,{0.42+0.055*\x+0.006*\x*\x+2*\hB})
		--
		plot[domain=\xR:\xL,samples=100,variable=\x]
			(\x,{0.42+0.055*\x+0.006*\x*\x+\hB/2})
		-- cycle;
	\draw[dashed,blue!70!black]
		plot[domain=\xL:\xR,samples=100,variable=\x]
			(\x,{0.42+0.055*\x+0.006*\x*\x+\hB/2});
	\draw[dashed,blue!70!black]
		plot[domain=\xL:\xR,samples=100,variable=\x]
			(\x,{0.42+0.055*\x+0.006*\x*\x+2*\hB});
	\foreach \xa/\xb in {0.720/1.501,1.501/2.192,2.192/2.819,2.819/3.397,3.397/3.935,3.935/4.442,4.442/4.900}{
		\pgfmathsetmacro{\pa}{0.42+0.055*\xa+0.006*\xa*\xa}
		\pgfmathsetmacro{\pb}{0.42+0.055*\xb+0.006*\xb*\xb}
		\draw[blue!65!black,line width=0.45pt,fill=blue!8,fill opacity=0.45]
			(\xa,{\pa+\hB/2}) rectangle (\xb,{\pb+2*\hB});
	}

	\path[fill=orange!12,fill opacity=0.55]
		plot[domain=\xL:\xR,samples=100,variable=\x]
			(\x,{0.42+0.055*\x+0.006*\x*\x+2*\hC})
		--
		plot[domain=\xR:\xL,samples=100,variable=\x]
			(\x,{0.42+0.055*\x+0.006*\x*\x+\hC/2})
		-- cycle;
	\draw[dashed,orange!80!black]
		plot[domain=\xL:\xR,samples=100,variable=\x]
			(\x,{0.42+0.055*\x+0.006*\x*\x+\hC/2});
	\draw[dashed,orange!80!black]
		plot[domain=\xL:\xR,samples=100,variable=\x]
			(\x,{0.42+0.055*\x+0.006*\x*\x+2*\hC});
	\foreach \xa/\xb in {0.720/1.124,1.124/1.501,1.501/1.856,1.856/2.192,2.192/2.512,2.512/2.819,2.819/3.113,3.113/3.397,3.397/3.670,3.670/3.935,3.935/4.192,4.192/4.442,4.442/4.685,4.685/4.900}{
		\pgfmathsetmacro{\pa}{0.42+0.055*\xa+0.006*\xa*\xa}
		\pgfmathsetmacro{\pb}{0.42+0.055*\xb+0.006*\xb*\xb}
		\draw[orange!80!black,line width=0.45pt,fill=orange!12,fill opacity=0.50]
			(\xa,{\pa+\hC/2}) rectangle (\xb,{\pb+2*\hC});
	}

	\draw[thick]
		plot[domain=\xL:\xR,samples=100,variable=\x]
			(\x,{0.42+0.055*\x+0.006*\x*\x})
		node[right] {$y=\phi_\ell(x)$};
	\node[purple!70!black,anchor=west] at
		(5.20,{0.42+0.055*5.20+0.006*5.20*5.20+1.40*\hA}) {$h_\mu$};
	\node[blue!70!black,anchor=west] at
		(5.20,{0.42+0.055*5.20+0.006*5.20*5.20+1.45*\hB}) {$h_{\mu+1}$};
	\node[orange!80!black,anchor=west] at
		(5.20,{0.42+0.055*5.20+0.006*5.20*5.20+1.20*\hC}) {$h_{\mu+2}$};

	\node at (3.15,-0.15) {\small (b) Coverings across dyadic transverse scales};
\end{scope}

\end{tikzpicture}
}%
\caption{
	Anisotropic Whitney-type covering by admissible boxes near the curve \(y=\phi_\ell(x)\).
	(a) A curved piece \(R_{\nu,h}\) at a transverse distance \(Y\approx h\) and its associated admissible box \(B_{\nu,h}\).
	(b) The admissible coverings at three consecutive dyadic transverse scales \(h_\mu\), \(h_{\mu+1}\), and \(h_{\mu+2}\). For each fixed transverse scale, the \(x\)- and \(y\)-projections of the boxes exhibit uniformly bounded overlap.
}
\label{fig:admissible_boxes}
\end{figure}

\Prop{Geometry of admissible boxes}{\label{prop:admissible_box}
For a fixed \(1 \le \ell \le K\), let \(\phi_\ell(x)=c_\ell x^\gamma\) with \(\gamma \neq 0\). Let \(I_{j_1} = [2^{-j_1-1}, \, 2^{-j_1+1}]\). Assume additionally that \(0<\epsilon_0 \le 2^{-5|\gamma|-5} c_1|\gamma|\), and fix a transverse scale \(0<h\le 2^{|\gamma|+1}\epsilon_0 2^{-\gamma j_1}\). Define
	\[
		R_h \defeq \set[\bigg]{(x,y) : x \in I_{j_1}, \ \frac{h}{2} \le y-\phi_\ell(x) \le 2h }.
	\]
Then there exists a partition of \(I_{j_1}\) into consecutive subintervals \([x_\nu, x_{\nu+1}]\) such that
	\begin{equation} \label{eq:partition_nu}
		\frac{h}{2^3} \le |\phi_\ell(x_{\nu+1}) - \phi_\ell(x_\nu)| \le \frac{h}{2^2}.
	\end{equation}
For each \(\nu\), define	the curved piece
	\[
		R_{\nu, h} \defeq 
		\set[\bigg]{
			(x,y): x \in [x_\nu, x_{\nu+1}], \ \frac{h}{2} \le y-\phi_\ell(x) \le 2h
			}
	\]
and the corresponding covering box
	\[
		B_{\nu,h} \defeq [x_\nu, x_{\nu+1}] \times 
			\bracket[\bigg]{
				\min_{x\in[x_\nu, x_{\nu+1}]} \phi_\ell(x)+\frac{h}{2}, \ 
				\max_{x\in[x_\nu, x_{\nu+1}]} \phi_\ell(x)+2h
			}.
	\]
We call \(B_{\nu,h}\) the \emph{admissible box} associated with \(R_{\nu,h}\).
Then the following hold.

\begin{enumerate}[label=\rm (A\arabic*)]
	\item \label{item:A1}
		The horizontal side length \(\Delta x_\nu\) of \(B_{\nu,h}\) satisfies
			\[
				\Delta x_\nu \approx 2^{-j_1(1-\gamma)}h.
			\]
		For every \((x,y) \in B_{\nu,h}\), write \(Y=y-\phi_\ell(x)\). Then the area of \(B_{\nu,h}\) satisfies
			\[
				|B_{\nu,h}| \approx x^{1-\gamma}Y^2.
			\]
	
	\item \label{item:A2}
		One has \(R_{\nu,h} \subset B_{\nu,h}\) for every \(\nu\), and hence \(R_h \subset \bigcup_\nu B_{\nu,h}\).
		The \(x\)-projections of the boxes \(\set{B_{\nu,h}}_\nu\) are at most \(2\)-overlapping, and their \(y\)-projections are at most \(14\)-overlapping.	
\end{enumerate}

Consequently, \(R_h\) is covered by axis-parallel boxes \(B_{\nu,h}\) whose \(x\)- and \(y\)-projections have uniformly bounded overlap and which satisfy
	\[
		|B_{\nu,h}| \approx |R_{\nu,h}| \quad \text{for every } \nu.
	\]
}
\pf{
It suffices to consider the case where \(\phi_\ell\) is increasing on \(I_{j_1}\); the decreasing case is identical. The assumed bounds on \(\epsilon_0\) and \(h\) guarantee that the total variation of \(\phi_\ell\) on \(I_{j_1}\) strictly exceeds \(h/2^3\). We therefore construct the partition satisfying \eqref{eq:partition_nu} by taking successive intervals over which the variation of \(\phi_\ell\) is exactly \(h/2^3\). If the variation on the final remaining interval is less than \(h/2^3\), we merge it with the preceding one.
Then the box \(B_{\nu,h}\) is
	\[
		B_{\nu,h} = [x_\nu, x_{\nu+1}] \times 
			\bracket[\bigg]{\phi_\ell(x_\nu)+\frac{h}{2}, \ \phi_\ell(x_{\nu+1})+2h}.
	\]

For \ref{item:A1}, since \(x\approx2^{-j_1}\) on \(I_{j_1}\), one has  \(|\phi_\ell^\prime(x)| \approx 2^{-j_1(\gamma-1)}\). Hence the mean value theorem and \eqref{eq:partition_nu} give
	\[
		h\approx 2^{-j_1(\gamma-1)} \Delta x_\nu.
	\]
Thus, we obtain the horizontal side length \(\Delta x_\nu \approx 2^{-j_1(1-\gamma)}h\).
Moreover, by \eqref{eq:partition_nu}, the vertical side length of \(B_{\nu,h}\) is comparable to \(h\), which yields
	\[
		|B_{\nu,h}| \approx \Delta x_\nu h \approx 2^{-j_1(1-\gamma)}h^2.
	\]
For \((x,y) \in B_{\nu,h}\), the definition of \(B_{\nu,h}\) and \eqref{eq:partition_nu} give
	\[
		\frac{h}{2^2} \le y-\phi_\ell(x) \le \phi_\ell(x_{\nu+1})-\phi_\ell(x_\nu) +2h \le \frac{h}{2^2}+2h.
	\]
Hence, setting \(Y=y-\phi_\ell(x)\), we obtain \(Y \approx h\). Since also \(x \approx 2^{-j_1}\), it follows that
	\[
		|B_{\nu,h}| \approx 2^{-j_1(1-\gamma)}h^2 \approx x^{1-\gamma}Y^2.
	\]

For \ref{item:A2}, the inclusion \(R_{\nu,h} \subset B_{\nu,h}\) follows directly from the definition of \(B_{\nu,h}\). Since the intervals \([x_\nu, x_{\nu+1}]\) partition \(I_{j_1}\), this gives \(R_h \subset \bigcup_\nu B_{\nu,h}\).

Moreover, the \(x\)-projections \([x_\nu, x_{\nu+1}]\) of \(B_{\nu,h}\) overlap only at their endpoints, so their overlap is at most \(2\).
 Let \(J_\nu = \bracket[\big]{\phi_\ell(x_\nu)+h/2, \ \phi_\ell(x_{\nu+1})+2h}\) be the \(y\)-projection of \(B_{\nu,h}\). If \(J_{\nu} \cap J_{\mu} \neq \emptyset\) for some \(\mu>\nu\), then the lower endpoint of \(J_\mu\) is no larger than the upper endpoint of \(J_\nu\), and hence \(\phi_\ell(x_\mu)+h/2 \le \phi_\ell(x_{\nu+1})+2h\). This gives \(\phi_\ell(x_\mu) - \phi_\ell(x_{\nu+1}) \le 3h/2\).
  Since \(\phi_\ell(x_{s+1})-\phi_\ell(x_s) \ge h/2^3\) by \eqref{eq:partition_nu}, this forces \((\mu-\nu-1)\frac{h}{2^3} \le \frac{3h}{2}\), which implies \(\mu-\nu \le 13\). Consequently, the \(y\)-projections are at most \(14\)-overlapping.
  
 Finally, by the definition of \(R_{\nu,h}\),
 	\[
 		|R_{\nu,h}| = \frac32 \Delta x_\nu h \approx |B_{\nu,h}|.
 	\]
 This completes the proof.
}
\Rmk{}{The choice of \(\epsilon_0\) in \eqref{eq:epsilon_separation}, together with the assumed upper bound on \(h\), ensures that every admissible box remains in \(\R_{>0}^2\) and does not intersect any other zero curve. As illustrated in \Cref{fig:admissible_boxes}, the associated admissible boxes \(B_{\nu,h}\) form an anisotropic Whitney-type covering of \(R_h\); their \(x\)- and \(y\)-projections have uniformly bounded overlap.
}

By \labelcref{item:A1} in \Cref{prop:admissible_box}, one has
	\[
		|B_{\nu,h}| \approx x^{1-\gamma}Y^2.
	\]
In view of the exponent vector \((1-\gamma,2)\) of \(x^{1-\gamma}Y^2\), this justifies the choice of \(\bv_\ell = (1-\gamma,2)\) in \eqref{eq:NP_near}.

\subparagraph{Horizontal curves.}
It remains to treat the case \(\gamma=0\). For a dyadic transverse scale \(0<h \lesssim 1\), set
	\[
		R_h \defeq \set[\bigg]{(x,y) : x \in I_{j_1}, \quad \frac h2 \le y-c_\ell \le 2h}.
	\]
Since the curve is horizontal, \(R_h\) is already an axis-parallel box in the original coordinates, and no further subdivision in \(x\) is needed. Writing \(Y=y-c_\ell\), for \((x,y) \in R_h\) one has \(x \approx 2^{-j_1}\) and \(Y \approx h\). Therefore
	\[
		|R_h| \approx 2^{-j_1}h.
	\]
In view of the exponent vector \((1,1)\) of \((2^{-j_1})^1h^1\), this justifies the choice of \(\bv_\ell=(1,1)\) in \eqref{eq:NP_near}.
This completes the derivation of the Newton pairs \((\P_\ell, \bv_\ell)\) declared in \labelcref{eq:NP_away,eq:NP_near}.


\subsection{Geometry of the Newton pairs} \label{subsec:geo_NP}

We now analyze the Newton polygons derived in \Cref{subsec:decomposition}.
The factorization
	\[
		P_{xy}(x,y) = u x^M y^N \prod_{\ell=1}^{L}(y-c_\ell x^\gamma)^{k_\ell}
	\]
makes the geometry of the Newton polygons essentially one-dimensional. Expanding this product shows that every exponent in \(\Lambda_0=\Lambda(P_{xy})\) lies on the line segment whose endpoints correspond to the lowest and highest powers of \(y\):
	\begin{equation} \label{eq:vertices_P0}
	\begin{alignedat}{2}
		\fm_{\textrm{bot}}
		&=(m_\mathrm{ bot,1},m_\mathrm{ bot,2})
		&&\defeq
		\paren[\Big]{
			M+\gamma\sum_{\ell=1}^L k_\ell,\ N
		},\\
		\fm_{\textrm{top}}
		&=(m_\mathrm{ top,1},m_\mathrm{ top,2})
		&&\defeq
		\paren[\Big]{
			M,\ N+\sum_{\ell=1}^L k_\ell
		}.
	\end{alignedat}
	\end{equation}
This segment is parallel to \((-\gamma,1)\).
	
Similarly, for \(1 \le \ell \le K\), set \(Y=y-c_\ell x^\gamma\). Recall that the exponent set \(\Lambda_\ell\) and the corresponding Newton polygon \(\P_\ell\) are both taken with respect to the adapted coordinates \((x,Y)\). Substituting \(y=Y+c_\ell x^\gamma\) into the factorization of \(P_{xy}\) interchanges the roles of \(N\) and \(k_\ell\). Hence, every exponent in \(\Lambda_\ell=\Lambda(P_{xy}(x,Y+c_\ell x^\gamma))\) lies on the line segment joining \(\fm_\ell\) and \(\fm_{\textrm{top}}\), where
	\begin{equation} \label{eq:vertices_Pl}
		\fm_{\ell} = (m_{\ell,1}, m_{\ell,2}) \defeq \paren[\Big]{M+\gamma N + \gamma \sum_{n \neq \ell} k_n, \ k_\ell}.
	\end{equation}
Here \(\fm_\ell\) is the endpoint corresponding to the lowest \(Y\)-degree in the adapted coordinates. This segment is also parallel to \((-\gamma,1)\).

Thus all the relevant exponent sets are collinear in the common direction \((-\gamma,1)\). The following lemma records the resulting geometry. To formulate the lemma, for \(\fq \in \R^2\setminus\set{\bm 0}\) and \(r \in \R\), we write
	\[
		\plane{\fq}{r} \defeq \set{\by \in \R^2 : \fq \cdot \by = r}.
	\]
 Recall that \(\conv E\) denotes the convex hull of a set \(E\).

\Lemma{Explicit description of the Newton polygons}{\label{lemma:NP_geometry}
Let \(\fm_{\mathrm{bot}}\), \(\fm_{\mathrm{top}}\), and \(\fm_\ell\) \((1 \le \ell \le K)\) be defined as in \labelcref{eq:vertices_P0,eq:vertices_Pl}. Then
	\[
		\P_0 = \conv\set{\fm_{\mathrm{bot}}, \fm_{\mathrm{top}}}, \qquad
		\P_\ell = \fm_\ell + \cone\set{(-\gamma,1)} \quad \text{ for } 1 \le \ell \le K.
	\]
Moreover, there exists \(r \in \R\) such that all these Newton polygons are contained in a common line \(\plane{\fq}{r}\), where \(\fq=(1,\gamma)\).
}
\pf{
Recall that \(\P_0 = \conv\paren{\Lambda_0}\) and \(\P_\ell = \conv\paren[\big]{\Lambda_\ell+\cone\set{(-\gamma,1)}}\).
Hence, by the preceding descriptions of \(\Lambda_0\) and \(\Lambda_\ell\), we have
	\[
		\P_0 = \conv\set{\fm_{\textrm{bot}}, \fm_{\textrm{top}}}, \qquad
		\P_\ell = \fm_\ell + \cone\set{(-\gamma,1)}.
	\]
For the second identity, we use \(\fm_{\textrm{top}} - \fm_\ell \in \cone\set{(-\gamma,1)}\).

Finally, set \(r = M+\gamma \paren[\big]{N+\sum_{n=1}^L k_n}\). Then for \(\fq=(1,\gamma)\), direct calculation gives
	\[
		\fq\cdot(-\gamma,1)=0, \qquad
		\fq\cdot \fm_{\textrm{bot}} = \fq \cdot \fm_{\textrm{top}} = \fq\cdot \fm_\ell
		= r.
	\]
Hence, all \(\P_\ell\) \((0 \le \ell \le K)\) are contained in \(\plane{\fq}{r}\).
}

\begin{figure}[htb]
\centering
\resizebox{0.72\textwidth}{!}{%
\begin{tikzpicture}[
	>=Latex,
	font=\small,
	axis/.style={
		->,
		line width=0.45pt,
		gray!70!black
	},
	polygon/.style={
		line width=1.8pt,
		black
	},
	coneray/.style={
		->,
		densely dashed,
		line width=0.9pt,
		black
	},
	areavector/.style={
		->,
		line width=1.2pt,
		black
	},
	dot/.style={
		circle,
		fill=black,
		inner sep=1.5pt
	},
	commonline/.style={
		line width=5pt,
		gray!30,
		line cap=round
	}
]


\begin{scope}[xshift=0cm,yshift=7cm]

\draw[axis]
	(-0.15,0) -- (5.15,0)
	node[right] {$m_1$};
\draw[axis]
	(0,-0.15) -- (0,5.35)
	node[above] {$m_2$};

\draw[commonline]
	(1.225,-0.10) -- (2.65,5.60)
	node[pos=.89,above left,text=gray!60!black] {$\plane{\fq}{r}$};

\draw[coneray]
	(0,0) -- (4.75,4.75)
	node[pos=.86,below right,xshift=-2pt] {$\cone(\bv_0)$};

\draw[areavector]
	(0,0) -- (0.75,0.75)
	node[pos=.80,below right] {$\bv_0$};

\coordinate (MB) at (1.375,0.50);
\coordinate (MT) at (2.375,4.50);

\draw[polygon]
	(MB) -- (MT)
	node[pos=.62,left] {$\P_0$};

\node[dot] at (MB) {};
\node[dot] at (MT) {};
\node[below right] at (MB) {$\fm_{\textrm{bot}}$};
\node[above right] at (MT) {$\fm_{\textrm{top}}$};

\node[dot] at (1.667,1.667) {};

\end{scope}

\begin{scope}[xshift=9cm,yshift=7cm]

\draw[axis]
	(-0.15,0) -- (5.15,0)
	node[right] {$m_1$};
\draw[axis]
	(0,-0.15) -- (0,5.35)
	node[above] {$m_2$};

\draw[commonline]
	(1.225,-0.10) -- (2.65,5.60)
	node[pos=.19,below right,text=gray!60!black] {$\plane{\fq}{r}$};

\draw[coneray]
	(0,0) -- (3.05,4.88)
	node[pos=.84,right,xshift=4pt] {$\cone(\bv_\ell)$};

\draw[areavector]
	(0,0) -- (0.625,1.00)
	node[pos=.75,right,xshift=2pt] {$\bv_\ell$};

\coordinate (ML) at (1.875,2.50);
\coordinate (RE) at (2.60,5.40);

\draw[polygon,->]
	(ML) -- (RE)
	node[pos=.48,left] {$\P_\ell$};

\node[dot] at (ML) {};
\node[below right] at (ML) {$\fm_\ell$};

\node[dot] at (2.083,3.333) {};

\end{scope}

\node at (6.45,8.80)
	{\Large $(\P_0, \bv_0)\;\Longrightarrow\;(\P_\ell,\bv_\ell)$};
\node at (7.12,6.20)
	{\footnotesize (a) Comparison between $(\P_0,\bv_0)$ and $(\P_\ell, \bv_\ell)$ for $\gamma\neq-1$};


\begin{scope}[xshift=0cm,yshift=0cm]

\draw[axis]
	(-0.15,0) -- (4.45,0)
	node[right] {$m_1$};
\draw[axis]
	(0,-0.15) -- (0,4.85)
	node[above] {$m_2$};

\draw[commonline]
	(-0.10,-0.10) -- (4.55,4.55);

\draw[coneray]
	(0,0) -- (4.50,4.50)
	node[pos=.99,above] {$\cone(\bv_0)$};

\draw[areavector]
	(0,0) -- (0.50,0.50)
	node[below right] {$\bv_0$};

\coordinate (MBm) at (1.05,1.05);
\coordinate (MTm) at (3.50,3.50);

\draw[polygon]
	(MBm) -- (MTm)
	node[pos=.58,above left] {$\P_0$};

\node[dot] at (MBm) {};
\node[dot] at (MTm) {};
\node[below right] at (MBm) {$\fm_{\textrm{bot}}$};
\node[below right] at (MTm) {$\fm_{\textrm{top}}$};

\end{scope}

\begin{scope}[xshift=9cm,yshift=0cm]

\draw[axis]
	(-0.15,0) -- (4.45,0)
	node[right] {$m_1$};
\draw[axis]
	(0,-0.15) -- (0,4.85)
	node[above] {$m_2$};

\draw[commonline]
	(-0.10,-0.10) -- (4.55,4.55);

\draw[coneray]
	(0,0) -- (4.50,4.50)
	node[pos=.99,above] {$\cone(\bv_\ell)$};

\draw[areavector]
	(0,0) -- (0.85,0.85)
	node[below right] {$\bv_\ell$};

\coordinate (MLm) at (1.75,1.75);
\coordinate (REm) at (4.10,4.10);

\draw[polygon,->]
	(MLm) -- (REm)
	node[pos=.55,above left] {$\P_\ell$};

\node[dot] at (MLm) {};
\node[below right] at (MLm) {$\fm_\ell$};

\end{scope}

\node at (6.45,1.8)
	{\Large $(\P_0,\bv_0)\;\Longrightarrow\;(\P_\ell,\bv_\ell)$};
\node at (7.12,-0.8)
	{\footnotesize (b) Comparison between $(\P_0,\bv_0)$ and $(\P_\ell,\bv_\ell)$ for $\gamma=-1$};

\end{tikzpicture}%
}
\caption{
Geometry of balanced Newton pairs.
(a) When \(\gamma \neq -1\), the rays \(\cone(\bv_0)\) and \(\cone(\bv_\ell)\) intersect the common line \(\plane{\fq}{r}\) transversely, so the forward and backward Newton distances of each pair coincide.
(b) When \(\gamma=-1\), the common line has the same direction as both \(\bv_0\) and \(\bv_\ell\). The polygon \(\P_0\) is a segment contained in \(\cone(\bv_0)\), whereas \(\P_\ell\) is an unbounded ray contained in \(\cone(\bv_\ell)\); hence \(\db(\P_\ell,\bv_\ell)=+\infty\).
}
\label{fig:newton_pairs}
\end{figure}
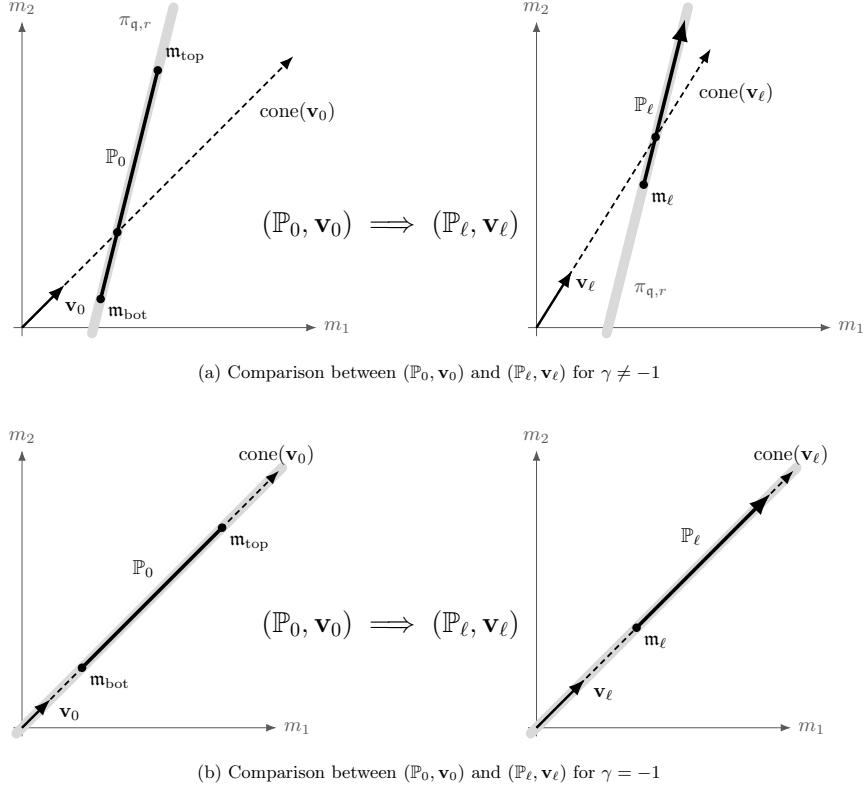

Recall from \labelcref{eq:NP_away,eq:NP_near} that
	\[
		\bv_0=(1,1), \qquad
		\bv_\ell =
		\begin{cases}
			(1-\gamma,2) \quad & \text{if } \gamma \neq 0, \\
			(1,1) & \text{if } \gamma =0,
		\end{cases}
		\quad \text{for } 1 \le \ell \le K.
	\]
Combining the explicit vectors above with \Cref{lemma:NP_geometry} immediately yields the following corollary; see \Cref{fig:newton_pairs}.

\Cor{}{\label{cor:ray_geometry}
The following statements hold.
\begin{enumerate}[\rm (1)]
	\item 
		For \(1\le \ell \le K\), if \((\P_\ell,\bv_\ell)\) is balanced, then \(\df(\P_\ell, \bv_\ell)>0\).
	\item
		If \(\gamma \neq -1\), then \(\df(\P_\ell,\bv_\ell)=\db(\P_\ell,\bv_\ell)\) for every balanced pair \((\P_\ell,\bv_\ell)\).
	\item
		If \(\gamma =-1\) and \((\P_0,\bv_0)\) is balanced, then every Newton pair is balanced and
			\[
				\P_\ell = [\df(\P_\ell, \bv_\ell), \db(\P_\ell, \bv_\ell)]\bv_\ell \quad \text{ for } 0 \le \ell \le K.
			\]
		Moreover, \(\db(\P_\ell, \bv_\ell)=+\infty\) for \(1 \le \ell \le K\).
\end{enumerate}
}

\subsection{Hessian bounds on admissible and testing boxes}

The geometry of Newton polygons in \Cref{lemma:NP_geometry} and the box scales in \Cref{prop:admissible_box} give the Hessian bounds used in the upper-bound estimates and in the testing arguments. We first record the bounds on the admissible covering boxes.

\Lemma{Hessian bounds on admissible boxes}{\label{lemma:admissible_box_hessian}
In the setting of \Cref{prop:admissible_box}, write \(Y=y-c_\ell x^\gamma\). Then on \(B_{\nu,h}\)
	\[
		|P_{xy}(x,y)| \approx |x|^{m_{\ell,1}}|Y|^{m_{\ell,2}},
		\qquad \text{and} \qquad
		|\partial_y^b P_{xy}(x,y)| \lesssim |x|^{m_{\ell,1}}|Y|^{m_{\ell,2}-b} \quad \text{for } b=1,2.
	\]
}
\pf{
By \Cref{prop:admissible_box}, one has \(x \approx 2^{-j_1}\) and \(|Y| \approx h\) on \(B_{\nu,h}\). The assumed upper bound on \(h\) therefore gives \(|Y| \lesssim \epsilon_0 x^\gamma\). The choice of \(\epsilon_0\) then gives
	\[
		y \approx x^\gamma, \quad \text{and} \quad
		|y-c_n x^\gamma| \approx x^\gamma
		\quad \text{for every } n \neq \ell.
	\]
Combining these estimates with the factorization
	\[
		P_{xy}(x,y)=ux^My^N \prod_{n=1}^L (y-c_n x^\gamma)^{k_n}
	\]
yields the normal-crossing estimate
	\[
		|P_{xy}(x,y)| 
		\approx |x|^{M+\gamma N + \gamma \sum_{n \neq \ell} k_n} |Y|^{k_\ell}
		= |x|^{m_{\ell,1}}|Y|^{m_{\ell,2}}.
	\]
Here we used \(\fm_\ell = (m_{\ell,1}, m_{\ell,2}) = (M+\gamma N + \gamma \sum_{n \neq \ell} k_n, k_\ell)\).

Finally, since \(\partial_y = \partial_Y\) with \(x\) fixed, differentiating the same factorization and using \(x^{-\gamma} \lesssim |Y|^{-1}\), which again follows from \(|Y| \lesssim \epsilon_0 x^\gamma\), yields
	\[
		|\partial_y^b P_{xy}(x,y)| \lesssim |x|^{m_{\ell,1}} |Y|^{m_{\ell,2}-b} \quad \text{for } b=1,2. 
		\qedhere
	\]
}

The preceding lemma gives the Hessian bounds on the admissible boxes used in the upper-bound argument. For the lower-bound and the unboundedness arguments, we use testing boxes at the corresponding side-length scales. Recall that the exponent set \(\Lambda_\ell\) is taken with respect to the adapted coordinates \((x,Y)\). We record the following elementary bounds for use in both \Cref{sec:lower_bounds,sec:unboundedness}. The additional hypothesis in \ref{item:T2} will be verified separately in each application.

\Lemma{Elementary bounds on testing boxes}{\label{lemma:box_hessian}
Fix \(0< \epsilon \ll 1\) sufficiently small.
\begin{enumerate}[label=\rm (T\arabic*)]
	\item \label{item:T1}
		Let \(x_0, y_0 > 0\) and define the box
			\[
				Q \defeq \set{(x,y):|x-x_0|\le \epsilon x_0, \ |y-y_0|\le \epsilon y_0}.
			\]
		Then
			\[
				\sup_{(x,y) \in Q} |P_{xy}(x,y)| \lesssim 
				\max_{(m_1, m_2) \in \Lambda_0} |x_0|^{m_1}|y_0|^{m_2}.
			\]
	\item \label{item:T2}
		Fix \(1 \le \ell \le K\).
		Let \(x_0, Y_0>0\) and set \(y_0 = c_\ell x_0^\gamma + Y_0\).
		Define the box
			\[
				Q \defeq \set{(x,y) : |x-x_0| \le \Delta x, \ |y-y_0| \le \Delta y},
			\]
		where
			\[
				\Delta x =
					\begin{cases}
						\epsilon x_0^{1-\gamma}Y_0, & \text{if } \gamma \neq 0, \\[0.5ex]
						\epsilon x_0, & \text{if } \gamma =0,
					\end{cases}
				\qquad
				\Delta y = \epsilon Y_0.
			\]
		If \(\gamma \neq 0\), assume additionally \(x_0^{-\gamma}Y_0 \le 1\). Then
			\[
				\sup_{(x,y) \in Q} |P_{xy}(x,y)| 
				\lesssim \max_{(m_1, m_2) \in \Lambda_\ell} |x_0|^{m_1}|Y_0|^{m_2}.
			\]
\end{enumerate}
}
\pf{
For \ref{item:T1}, one has \(x \approx x_0\) and \(y \approx y_0\) on \(Q\). The stated estimate therefore follows immediately by expanding \(P_{xy}\) into monomials.

For \ref{item:T2}, set \(Y=y-c_\ell x^\gamma\). 
The construction of \(Q\) and the mean value theorem yield \(x \approx x_0\) and \(Y \approx Y_0\) on \(Q\). Thus, using the expansion
	\(
		P_{xy}(x,Y+c_\ell x^\gamma) 
		= \sum_{(m_1, m_2) \in \Lambda_\ell} \widetilde c_{m_1,m_2} x^{m_1}Y^{m_2}
	\)
in the adapted coordinates \((x,Y)\), we obtain
	\[
		\sup_{(x,y)\in Q}|P_{xy}(x,y)|
		\lesssim
		\max_{(m_1,m_2)\in\Lambda_\ell}
		|x_0|^{m_1}|Y_0|^{m_2}.
		\qedhere
	\]
}


\subsection{Analytic lemmas}

We record two operator estimates that will be applied to dyadic blocks and admissible boxes. The first uses only the size of the kernel support, while the second uses oscillation through a lower bound for the Hessian.

We use the following size estimate; see \cite[Lemma~1.6]{PS94}.

\Lemma{Size estimate, \cite{PS94}}{\label{lemma:size_est}
Let \(T\) be an integral operator with kernel \(K(x,y)\), and assume that
	\begin{enumerate}[label=\rm (\alph*)]
		\item \(|K(x,y)| \le 1\),
		\item for each \(y\), the set \(\set{x: K(x,y) \neq 0}\) has measure at most \(\Delta_x\),
		\item for each \(x\), the set \(\set{y: K(x,y) \neq 0}\) has measure at most \(\Delta_y\).
	\end{enumerate}
Then \(\norm{T}_{\op} \le (\Delta_x\Delta_y)^{1/2}\).
}

This follows from Schur's test; we omit the proof.
We also use the following oscillatory estimate, which is a variant of the operator van der Corput lemma of Phong--Stein \cite{PS97}. Its proof follows that of \cite[Lemma~1.1]{PS94}, replacing the polynomiality and Hessian bounds with our assumptions along horizontal segments. We include the details to keep the exposition self-contained.

\Lemma{Oscillatory estimate}{\label{lemma:osc_est}
Let \(T^\lambda\) be an oscillatory integral operator of the form
	\[
		T^\lambda f(x) = \int \ee^{\ii\lambda P(x,y)} \psi_R(x,y) f(y) \dy.
	\]
Assume that
	\begin{enumerate}[label=\rm (\alph*)]
		\item \(\psi_R(x,y)\) is supported in a box \(R\) of size \(\Delta_x \times \Delta_y\),
		\item \(|\partial_y^b \psi_R(x,y)| \le C(\Delta_y)^{-b}\) for \(b=0,1,2\),
		\item whenever \((x,y), (x^\prime,y) \in \supp \psi_R\), \(|P_{xy}(z,y)| \ge \omega >0\) for every \(z\) between \(x\) and \(x^\prime\),
		\item under the same conditions, \(|\partial_y^b P_{xy}(z,y)| \le C\omega (\Delta_y)^{-b}\) for \(b=1,2\).
	\end{enumerate}
Then \(\norm{T^\lambda}_{\op} \lesssim (\lambda \omega)^{-1/2}\) with a constant depending only on \(C\).
}
\pf{
The operator \(T^\lambda(T^\lambda)^*\) is an integral operator with kernel \(K(x,x^\prime)\) given by
	\[
		K(x,x^\prime) 
		= \int \ee^{\ii\lambda\Phi(x,x^\prime,y)}\psi_R(x,y)\overline{\psi_R(x^\prime,y)}\dy,
		\quad \text{where } \Phi(x,x^\prime,y) = P(x,y)-P(x^\prime,y).
	\]
We establish the following estimates for \(K(x,x^\prime)\):
	\begin{subequations} \label{eq:subeq_Osc}
		\begin{align}
			|K(x,x^\prime)| &\le C\Delta_y, 
				\label{eq:subeq_Osc1} \\
			|K(x,x^\prime)| &\le C\Delta_y^{-1}(\lambda \omega |x-x^\prime|)^{-2}. \label{eq:subeq_Osc2} 	
		\end{align}
	\end{subequations}
The estimate \eqref{eq:subeq_Osc1} is a consequence of (a) and (b).

To obtain \eqref{eq:subeq_Osc2}, we use integration by parts. By (c) and continuity, \(P_{xy}\) does not change sign in the interval \([x^\prime, x]\), and thus
	\begin{equation} \label{eq:LB_derivative_Osc}
		|\partial_y\Phi(x,x^\prime,y)| 
		= \abs[\bigg]{\int_{x^\prime}^x P_{xy}(z,y) \d{z}}
		\ge \omega|x-x^\prime|.
	\end{equation}
Moreover, by the assumption (d), for \(b=1,2\),
	\begin{equation} \label{eq:UB_derivative_Osc}
		|\partial_y^{b+1} \Phi(x,x^\prime,y)| \le C\omega|x-x^\prime|(\Delta_y)^{-b}.
	\end{equation}
Using the operator \(\cL=\frac{1}{\ii\lambda\partial_y\Phi}\partial_y\), which satisfies \(\cL \ee^{\ii\lambda\Phi}=\ee^{\ii\lambda\Phi}\), we rewrite the kernel as
	\[
		K(x,x^\prime)
		= \int \ee^{\ii\lambda\Phi(x,x^\prime,y)} 
			(\cL^\intercal)^2\paren[\big]{\psi_R(x,y)\overline{\psi_R(x^\prime,y)}} \dy,
		\quad \text{where } 
		\cL^\intercal f = \partial_y\paren[\Big]{\frac{-f}{\ii\lambda\partial_y\Phi}}.
	\]
Setting \(A(x,x^\prime,y)=\psi_R(x,y)\overline{\psi_R(x^\prime,y)}\) and integrating by parts twice in \(y\) gives
	\[
		|K(x,x^\prime)|
		\lesssim \lambda^{-2} \int \paren[\bigg]{
			\frac{|\partial_y^2 A|}{|\partial_y \Phi|^2}
			+ \frac{|\partial_y A||\partial_y^2 \Phi|}{|\partial_y\Phi|^3}
			+ \frac{|A||\partial_y^3\Phi|}{|\partial_y\Phi|^3}
			+ \frac{|A||\partial_y^2\Phi|^2}{|\partial_y\Phi|^4}
		}\dy.
	\]
Combining \labelcref{eq:LB_derivative_Osc,eq:UB_derivative_Osc} with the assumption (b) yields \eqref{eq:subeq_Osc2}.

By \eqref{eq:subeq_Osc}, we have
	\[
		|K(x,x^\prime)| \le \frac{C\Delta_y}{1+(\lambda\omega\Delta_y|x-x^\prime|)^2},
	\]
which implies
	\[
		\sup_x \int |K(x,x^\prime)| \d{x^\prime} \le C(\lambda\omega)^{-1}, \qquad
		\sup_{x^\prime} \int |K(x,x^\prime)| \dx \le C(\lambda\omega)^{-1}.
	\]
Finally, by Schur's test, we obtain \(\norm{T^\lambda(T^\lambda)^*}_{\op} \lesssim (\lambda\omega)^{-1}\), which proves the result.
}


\section{Upper bounds for balanced Newton pairs} \label{sec:upper_bounds}

Throughout this section, we assume that the family of Newton pairs \(\fN(P)\), defined in \eqref{eq:NF}, is balanced.
For \(0 \le \ell \le K\), we define
	\begin{equation} \label{eq:delta_l}
		\delta_\ell(\lambda) \defeq
		\begin{cases}
			\df(\P_\ell,\bv_\ell), & \text{if } \lambda \ge 1, \\[0.5ex]
			\db(\P_\ell,\bv_\ell), & \text{if } 0< \lambda <1.
		\end{cases}
	\end{equation}
When \(\delta_\ell(\lambda)=+\infty\), the quantity \(1/(\delta_\ell(\lambda)+1)\) is understood as \(0\).	Recall from \Cref{subsec:decomposition} that
	\[
		T^\lambda = T^\lambda_0 + \sum_{\ell=1}^K T^\lambda_\ell.
	 \]
In this section, we estimate each piece \(T^\lambda_\ell\) in terms of \(\delta_\ell(\lambda)\). We first estimate \(T^\lambda_0\) and then the pieces \(T^\lambda_\ell\) for \(1 \le \ell \le K\).


\subsection{Away from the zero curves} \label{subsec:UB_away}

Since the Newton pair \((\P_0, \bv_0)\) is balanced, the definition in \eqref{eq:newton_distance} yields
	\[
		\P_0 \cap \cone(\bv_0) = \bracket[\big]{\df(\P_0,\bv_0), \db(\P_0,\bv_0)}\bv_0.
	\]
We establish
	\[
		\norm{T^\lambda_0}_{\op} \lesssim \lambda^{-\frac{1}{2(\delta_0(\lambda)+1)}},
	\]
where \(\delta_0(\lambda)\) is the Newton distance defined in \eqref{eq:delta_l}.

\subsubsection{Away geometry and dyadic decomposition}

We first set up the geometry and dyadic decomposition needed for the block estimates and the Cotlar--Stein summation. By \eqref{eq:supp_fts},
	\[
		\supp\, \psi_0 \subset 
		\bigcap_{\ell=1}^K \set[\Big]{ (x,y) : |y-c_\ell x^\gamma| \ge \frac{\epsilon_0}{2} x^\gamma}.
	\]
Enlarging this separated region slightly, we set
	\begin{equation} \label{eq:Omega0@away}
		\Omega_0 
		\defeq \bigcap_{\ell=1}^K \set[\Big]{(x,y): |y-c_\ell x^{\gamma}| \ge \frac{\epsilon_0}{4}x^{\gamma}},
	\end{equation}
so that \(\supp\, \psi_0 \subset \Omega_0\). When \(K=0\), we set \(\Omega_0=\R_{>0}^2\). On \(\Omega_0\), one has \(|y-c_nx^\gamma| \approx |y|+|x|^\gamma\) for all \(1\le n\le L\). Substituting this into the factorization \eqref{eq:quasi_form} and using the vertices of \(\P_0\) in \eqref{eq:vertices_P0} reveals that the Hessian \(P_{xy}\) is normal-crossing:
	\begin{equation} \label{eq:mono_hess@away}
		|P_{xy}(x,y)| \approx 
			|x|^{m_\mathrm{ bot,1}}|y|^{m_\mathrm{ bot,2}} 
			+ |x|^{m_\mathrm{ top,1}}|y|^{m_\mathrm{ top,2}} 
			\quad \text{on } \Omega_0.
	\end{equation}
Moreover, for \(a,b \in \Z_{\ge0}\), differentiating the same factorization yields the symbolic bounds
	\begin{equation} \label{eq:symbolic_bounds@away}
		|\partial_x^a \partial_y^b P_{xy}(x,y)| \lesssim |x|^{-a}|y|^{-b}|P_{xy}(x,y)| 
		\quad \text{on } \Omega_0.
	\end{equation}
	
The Cotlar--Stein summation relies on off-diagonal estimates obtained via integration by parts along vertical and horizontal segments. To ensure that the Hessian \(P_{xy}\) does not change sign on these segments, we decompose \(\psi_0\) according to the regions between the curves:
	\[
		\psi_0=\sum_{s=0}^K \psi_{0,s}.
	\]
If \(K=0\), this subdivision is unnecessary, and we take \(\psi_{0,0}=\psi_0\) with \(U_0 \defeq \R_{>0}^2\). When \(K \ge 1\), we choose smooth cutoffs \(\psi_{0,s}\) depending only on \(yx^{-\gamma}\) such that \(\psi_{0,s}\) is supported in \(U_s\) and \(\psi_0 =\sum_{s=0}^K \psi_{0,s}\), where
	\[
		\begin{alignedat}{3}
			&U_0 &&\defeq \set[\Big]{(x,y) \in \R_{>0}^2 : 0<y<\paren[\Big]{c_1-\frac{\epsilon_0}{4}}x^\gamma}, \\[0.5ex]
			&U_s &&\defeq \set[\Big]{(x,y) \in \R_{>0}^2:
			\paren[\Big]{c_s+\frac{\epsilon_0}{4}}x^\gamma<y<\paren[\Big]{c_{s+1}-\frac{\epsilon_0}{4}}x^\gamma
			} \quad \text{for } 1 \le s \le K-1, \\[0.5ex]
			&U_K &&\defeq \set[\Big]{(x,y) \in \R_{>0}^2: y>\paren[\Big]{c_K+\frac{\epsilon_0}{4}}x^\gamma}.
		\end{alignedat}
	\]
By construction, each \(U_s\) is contained in a connected component of \(\R_{>0}^2 \setminus \bigcup_{\ell=1}^K \set{y=c_\ell x^\gamma}\), with a uniform separation from the boundary curves.

The following elementary geometric observation ensures that the Hessian \(P_{xy}\) does not change sign on the vertical and horizontal segments used below.

\Lemma{}{\label{lemma:curve_separated}
For each \(0 \le s \le K\), let \(U_s\) be the open region defined above. Then the following hold.
\begin{enumerate}[label=\rm (\arabic*)]
	\item 
		If \((x,y)\) and \((x,y^\prime)\) belong to \(U_s\), then the vertical segment connecting them lies in \(U_s\).
	\item
		If \((x,y)\) and \((x^\prime,y)\) belong to \(U_s\), then the horizontal segment connecting them lies in \(U_s\).
\end{enumerate}
Consequently, \(P_{xy}\) has a constant sign on these segments.
}
\pf{
When \(K=0\), the assertions (1) and (2) are trivial since \(U_0 = \R_{>0}^2\); we thus assume \(K \ge 1\) in proving them.

For fixed \(x>0\), the defining inequalities of \(U_s\) restrict \(y\) to an interval, proving (1). For fixed \(y>0\), they restrict \(x^\gamma\) to an interval. If \(\gamma=0\), this imposes no constraint on \(x\) and the horizontal sections are intervals; otherwise \(x \mapsto x^\gamma\) is strictly monotone on \(\R_{>0}\), so the corresponding \(x\)-values again form an interval.

Finally, since \eqref{eq:Omega0@away} gives \(U_s \subset \Omega_0\), the region \(U_s\) contains no real zeros of \(P_{xy}\). By continuity, \(P_{xy}\) maintains a constant sign on any segment contained in \(U_s\).
}

Corresponding to this spatial decomposition, we define
	\[
		T^\lambda_{0,s} f(x) = \int \ee^{\ii\lambda P(x,y)} \psi_{0,s}(x,y) f(y) \dy 
		\quad \text{for } 0 \le s \le K,
	\]
so that \(T^\lambda_0 = \sum_{s=0}^K T^\lambda_{0,s}\). It suffices to estimate \(T^\lambda_{0,s}\) uniformly in \(s\). For a fixed \(s\), we decompose
	\[
		T^\lambda_{0,s} = \sum_{J \in \Z^2} T^\lambda_{J,s} \quad \text{for } J=(j_1, j_2),
	\]
where
	\[
		T^\lambda_{J,s} f(x)
		= \int \ee^{\ii\lambda P(x,y)}\psi_{0,s}(x,y)
			\chi\paren[\Big]{\frac{x}{2^{-j_1}}}\chi\paren[\Big]{\frac{y}{2^{-j_2}}}f(y)\dy.
	\]
We fix \(s\) for the remainder of \Cref{subsec:UB_away} and write \(T^\lambda_J = T^\lambda_{J,s}\).

For a vertex \(\fm\) of a polygon \(\P\), its \emph{dual face}\footnote{This is also called the normal cone of \(\P\) at \(\fm\).} \(\fm^\vee\) is defined as the cone of inward normal vectors to the supporting lines of \(\P\) through \(\fm\):
	\[
		\fm^\vee \defeq \set{\fq \in \R^2: \fq \cdot \fm = \min_{\fp \in \P} \fq \cdot \fp}.
	\]
Since \(\P_0=\conv\set{\fm_{\textrm{bot}}, \fm_{\textrm{top}}}\), these dual faces are given by
	\[
		\fm_{\textrm{bot}}^\vee
			= \set[\big]{\fq \in \R^2: \fq \cdot \fm_{\textrm{bot}} \le \fq \cdot \fm_{\textrm{top}}} 
		\quad \text{and} \quad
		\fm_{\textrm{top}}^\vee
			= \set[\big]{\fq \in \R^2: \fq \cdot \fm_{\textrm{top}} \le \fq \cdot \fm_{\textrm{bot}}},
	\]
so that \(\R^2= \fm_{\textrm{bot}}^\vee \cup \fm_{\textrm{top}}^\vee\). Assigning their intersection to \(\fm_{\textrm{bot}}^\vee\), we obtain the disjoint decomposition
	\[
		\R^2 
		= \fm_{\textrm{bot}}^\vee \sqcup 
		\paren[\big]{ \fm_{\textrm{top}}^\vee \setminus \fm_{\textrm{bot}}^\vee}.
	\]
Accordingly, we decompose the operator \(T^\lambda_{0,s}\) as
	\[
		T^\lambda_{0,s} 
		= T^\lambda_\mathrm{ bot} + T^\lambda_\mathrm{ top}
		\defeq
		\sum_{J \in \Z^2 \cap \fm_{\textrm{bot}}^\vee} T^\lambda_J 
			 + \sum_{J \in \Z^2 \cap (\fm_{\textrm{top}}^\vee \setminus \fm_{\textrm{bot}}^\vee)} T^\lambda_J.
	\]
It is enough to estimate \(T^\lambda_\mathrm{ bot}\); the estimate for \(T^\lambda_\mathrm{ top}\) follows symmetrically, with \(\fm_{\textrm{top}}\) replacing \(\fm_{\textrm{bot}}\).

\subsubsection{Block estimates}

We now prove the uniform block estimate needed for the Cotlar--Stein summation:
	\begin{equation} \label{eq:uniform_block@away}
		\norm{T^\lambda_J}_{\op}	\lesssim
		\lambda^{-\frac{1}{2(\delta_0(\lambda)+1)}},
		\quad \text{where } J\in\Z^2 \cap \fm_{\textrm{bot}}^\vee.
	\end{equation}
Let \(J \in \Z^2 \cap \fm_{\textrm{bot}}^\vee\). Then the kernel of \(T^\lambda_J\) is supported in a box where \(x \approx 2^{-j_1}\) and \(y\approx 2^{-j_2}\). Since \(\bv_0=(1,1)\) by \eqref{eq:NP_away}, this box has area \(\approx 2^{-J\cdot \bv_0}\), and \Cref{lemma:size_est} yields the size estimate
	\[
		\norm[\big]{T^\lambda_J}_{\op} \lesssim 2^{-J\cdot \bv_0/2}.
	\]

For the oscillatory estimate, we verify the hypotheses of \Cref{lemma:osc_est} along horizontal segments. Let \((x,y)\) and \((x^\prime,y)\) belong to the support of the kernel of \(T^\lambda_J\). By \Cref{lemma:curve_separated}, the horizontal segment joining these points lies in \(U_s \subset \Omega_0\). If \(z\) lies on this segment, then \(z \approx 2^{-j_1}\) and \(y \approx 2^{-j_2}\). Since \(J \in \Z^2 \cap \fm_{\textrm{bot}}^\vee\), one has \(J\cdot \fm_{\textrm{bot}} \le J \cdot \fm_{\textrm{top}}\), and therefore \eqref{eq:mono_hess@away} gives
	\[
		|P_{xy}(z,y)| \approx 2^{-J\cdot \fm_{\textrm{bot}}}.
	\]
Combined with the symbolic bounds \eqref{eq:symbolic_bounds@away}, this gives
	\[
		|\partial_y^b P_{xy}(z,y)| \lesssim 2^{-J\cdot \fm_{\textrm{bot}}}(2^{-j_2})^{-b} \quad \text{for } b=1,2.
	\]
The cutoff factors satisfy the same \(y\)-scale bounds, since the derivatives of \(\psi_{0,s}\) are supported where \(yx^{-\gamma} \approx 1\). Therefore, \Cref{lemma:osc_est} with \(\omega=2^{-J\cdot\fm_{\textrm{bot}}}\) yields
	\[
		\norm[\big]{T^\lambda_J}_{\op} \lesssim (\lambda 2^{-J\cdot \fm_{\textrm{bot}}})^{-1/2}.
	\]
Combining these two estimates, we obtain
	\begin{equation} \label{eq:block_norm@away}
		\norm[\big]{T^\lambda_J}_{\op} \lesssim 
		 \min\set[\Big]{ 2^{-J\cdot \bv_0 /2} , (\lambda 2^{-J\cdot \fm_{\textrm{bot}}})^{-1/2}}. 
	\end{equation}

It remains to relate \(J\cdot \fm_{\textrm{bot}}\) to the Newton distance. By the definition of the dual face, \(J\in \fm_{\textrm{bot}}^\vee\) implies \(J \cdot \fm_{\textrm{bot}} = \min_{\fp \in \P_0} J \cdot \fp\). Since \(\delta_0(\lambda)\bv_0 \in \P_0\), it follows that \(J\cdot \fm_{\textrm{bot}} \le \delta_0(\lambda) J\cdot \bv_0\).
Substituting this into \eqref{eq:block_norm@away} and balancing the two quantities in the minimum yields
	\begin{equation} \label{eq:bound_balance@away}
		\norm[\big]{T^\lambda_J}_{\op} \lesssim
		\min\set[\Big]{ 2^{-J\cdot \bv_0/2}, \lambda^{-1/2}2^{\delta_0(\lambda)J\cdot\bv_0/2}}
		\lesssim \lambda^{-\frac{1}{2(\delta_0(\lambda)+1)}},
	\end{equation}
which is \eqref{eq:uniform_block@away}.

\subsubsection{Cotlar--Stein summation}

To sum the family \(\set[\big]{T^\lambda_J}_{J \in \Z^2 \cap \fm_{\textrm{bot}}^\vee}\), we establish the Cotlar--Stein estimates
	\begin{subequations} \label{eq:CS_claims@away}
		\begin{align}
			\norm[\big]{(T^\lambda_J)^*T^\lambda_{J^\prime}}_{\op}
				&\lesssim \lambda^{-\frac{1}{\delta_0(\lambda)+1}}2^{-\frac{1}{2}|J-J^\prime|}, 
				\label{eq:CS_claim1@away} \\
			\norm[\big]{T^\lambda_J(T^\lambda_{J^\prime})^*}_{\op}
				&\lesssim \lambda^{-\frac{1}{\delta_0(\lambda)+1}}2^{-\frac{1}{2}|J-J^\prime|}
				\label{eq:CS_claim2@away}
		\end{align}
	\end{subequations}
for all \(J, J^\prime \in \Z^2 \cap \fm_{\textrm{bot}}^\vee\), where \(|J-J^\prime| = |j_1-j_1^\prime|+|j_2-j_2^\prime|\). Once \eqref{eq:CS_claims@away} is proved, the Cotlar--Stein lemma \cite{St93} gives
	\[
		\norm[\Big]{\sum_{J\in\Z^2 \cap \fm_{\textrm{bot}}^\vee}T^\lambda_J}_{\op}
		\lesssim
		\lambda^{-\frac{1}{2(\delta_0(\lambda)+1)}}.
	\]
	
For \((T^\lambda_J)^*T^\lambda_{J^\prime}\), the \(x\)-supports are disjoint unless \(|j_1-j_1^\prime| \le 2\), while for \(T^\lambda_J(T^\lambda_{J^\prime})^*\), the \(y\)-supports are disjoint unless \(|j_2-j_2^\prime|\le 2\). We isolate the off-diagonal estimates needed for \eqref{eq:CS_claims@away} in the following lemma.

\Lemma{Off-diagonal block estimates}{\label{lemma:off_diagonal@away}
Let \(J,J^\prime \in \Z^2 \cap \fm_{\textrm{bot}}^\vee\).
	\begin{enumerate}[\rm (1)]
		\item 
			If \(|j_1-j_1^\prime| \le 2\) and \(j_2^\prime \ge j_2 +3\), then
				\begin{equation} \label{eq:off_diagonal1@away}
					\norm[\big]{(T^\lambda_J)^*T^\lambda_{J^\prime}}_{\op}
					\lesssim
					2^{-\frac12|j_2-j_2^\prime|}
					\min\set[\Big]{
						2^{-J\cdot\bv_0}, (\lambda 2^{-J\cdot \fm_{\textrm{bot}}})^{-1}
					}.
				\end{equation}
		\item
			If \(|j_2-j_2^\prime| \le 2\) and \(j_1^\prime \ge j_1 +3\), then
				\begin{equation} \label{eq:off_diagonal2@away}
					\norm[\big]{T^\lambda_J(T^\lambda_{J^\prime})^*}_{\op}
					\lesssim
					2^{-\frac12|j_1-j_1^\prime|}
					\min\set[\Big]{
						2^{-J\cdot\bv_0}, (\lambda 2^{-J\cdot \fm_{\textrm{bot}}})^{-1}
					}.
				\end{equation}
	\end{enumerate}
}
\pf{
We prove \eqref{eq:off_diagonal1@away}; the same argument proves \eqref{eq:off_diagonal2@away}, with the roles of \(x\) and \(y\) interchanged.

The kernel of \((T^\lambda_J)^*T^\lambda_{J^\prime}\) is
	\[
		K(y,y^\prime) 
		= \int \ee^{-\ii\lambda \Phi(x,y,y^\prime)} A(x,y,y^\prime)\dx,
		\quad \text{where } \Phi(x,y,y^\prime) = P(x,y)-P(x,y^\prime).
	\]
Here, \(A(x,y,y^\prime)\) collects the spatial cutoffs, and
	\[
		|A(x,y,y^\prime)| \lesssim 1, \qquad |\partial_x A(x,y,y^\prime)| \lesssim 2^{j_1}.
	\]
On the support of \(A\), the vertical segment joining \((x,y^\prime)\) and \((x,y)\) lies in \(U_s\) by \Cref{lemma:curve_separated}. Hence \(P_{xy}\) has a constant sign on this segment. 

Since \(j_2^\prime \ge j_2+3\), the dyadic support conditions give \([y/2,\,y] \subset [y^\prime,\,y]\). For \(z \in [y/2,\,y]\), one has \(x\approx 2^{-j_1}\) and \(z \approx 2^{-j_2}\). Since \(J \in \Z^2 \cap \fm_{\textrm{bot}}^\vee\), one has \(J\cdot \fm_{\textrm{bot}} \le J \cdot \fm_{\textrm{top}}\), and hence  \eqref{eq:mono_hess@away} gives \(|P_{xy}(x,z)| \approx 2^{-J\cdot \fm_{\textrm{bot}}}\). The constant sign of \(P_{xy}\) then yields
	\[
		\abs{\partial_x\Phi(x,y,y^\prime)}
		= \abs[\bigg]{\int_{y^\prime}^y P_{xy}(x,z)\d{z}} 
		= \int_{y^\prime}^y \abs{P_{xy}(x, z)} \d{z} 
		\ge \int_{y/2}^y \abs{P_{xy}(x, z)} \d{z} 
		\approx 2^{-j_2}2^{-J\cdot\fm_{\textrm{bot}}}.
	\]
Moreover, since \(|z| \le |y|\) for \(z \in [y^\prime, y]\), the normal-crossing estimate \eqref{eq:mono_hess@away} yields
	\[
		\begin{alignedat}{2}
			|P_{xy}(x,z)| 
			&\lesssim |x|^{m_\mathrm{ bot,1}}|z|^{m_\mathrm{ bot,2}}+|x|^{m_\mathrm{ top,1}}|z|^{m_\mathrm{ top,2}} \\ 
			&\le |x|^{m_\mathrm{ bot,1}}|y|^{m_\mathrm{ bot,2}}+|x|^{m_\mathrm{ top,1}}|y|^{m_\mathrm{ top,2}} 
			\lesssim 2^{-J\cdot \fm_{\textrm{bot}}}.
		\end{alignedat}
	\]
Hence \eqref{eq:symbolic_bounds@away} yields
	\[
		\abs{\partial_x^2 \Phi (x,y,y^\prime)}
		\le \int_{y^\prime}^y |\partial_x P_{xy}(x,z)| \d{z}
		\lesssim |x|^{-1} \int_{y^\prime}^y \abs{P_{xy}(x,z)} \d{z}
		\lesssim 2^{j_1}2^{-j_2}2^{-J\cdot \fm_{\textrm{bot}}}.
	\]
Integrating by parts once in \(x\), using \(\cL=\frac{-1}{\ii\lambda\partial_x\Phi}\partial_x\), yields
	\[
		\begin{alignedat}{3}
			|K(y,y^\prime)|
			&= \abs[\bigg]{ 
					\int \ee^{-\ii\lambda \Phi} 
						\bracket[\bigg]{ 
							\partial_x \paren[\Big]{\frac{A}{\ii\lambda \partial_x \Phi}}
						}\dx 
					} \\[0.5ex]
			&\lesssim
			\lambda^{-1}
				\int \paren[\bigg]{
					\frac{|\partial_x A|}{|\partial_x\Phi|} + \frac{|A||\partial_x^2\Phi|}{|\partial_x\Phi|^2}
				}\dx  \\[0.5ex]
			&\lesssim
			\lambda^{-1} 2^{-j_1}
				\paren[\bigg]{
					\frac{2^{j_1}}{2^{-j_2}2^{-J\cdot\fm_{\textrm{bot}}}} 
					+ \frac{2^{j_1}2^{-j_2}2^{-J\cdot\fm_{\textrm{bot}}}} {(2^{-j_2}2^{-J\cdot\fm_{\textrm{bot}}})^2}
				}
			\lesssim 2^{j_2}(\lambda2^{-J\cdot\fm_{\textrm{bot}}})^{-1}.
		\end{alignedat}
	\]
Since \(y \approx 2^{-j_2}\) and \(y^\prime \approx 2^{-j_2^\prime}\) on the kernel support, Schur's test gives
	\[
		\norm[\big]{(T^\lambda_J)^*T^\lambda_{J^\prime}}_{\op}
		\lesssim 2^{-\frac{1}{2}|j_2-j_2^\prime|}(\lambda 2^{-J\cdot \fm_{\textrm{bot}}})^{-1}.
	\]
The pointwise size bound \(|K(y,y^\prime)| \lesssim 2^{-j_1}\) similarly gives
	\[
		\norm[\big]{(T^\lambda_J)^*T^\lambda_{J^\prime}}_{\op}
		\lesssim 2^{-\frac{1}{2}|j_2-j_2^\prime|}2^{-J\cdot\bv_0}.
	\]
Combining these estimates proves \eqref{eq:off_diagonal1@away}.
}

We now prove \eqref{eq:CS_claim1@away}. By the disjointness of the \(x\)-supports, we may assume \(|j_1-j_1^\prime|\le2\). When also \(|j_2-j_2^\prime|\le2\), the uniform block estimate \eqref{eq:uniform_block@away} is sufficient since \(|J-J^\prime|\le4\). Otherwise, by symmetry, we may assume \(j_2^\prime\ge j_2+3\). Then \eqref{eq:off_diagonal1@away} and the balancing argument as in \eqref{eq:bound_balance@away} give
		\[
		\norm[\big]{(T^\lambda_J)^*T^\lambda_{J^\prime}}_{\op}
		\lesssim \lambda^{-\frac{1}{\delta_0(\lambda)+1}} 2^{-\frac{1}{2}|j_2-j_2^\prime|}.
	\]
Since \(|j_1-j_1^\prime| \le 2\), the bound \eqref{eq:CS_claim1@away} follows.

We next prove \eqref{eq:CS_claim2@away}. By the disjointness of the \(y\)-supports, we may assume \(|j_2-j_2^\prime|\le2\). When also \(|j_1-j_1^\prime|\le2\), the uniform block estimate \eqref{eq:uniform_block@away} is again sufficient. Otherwise, by symmetry, we may assume \(j_1^\prime\ge j_1+3\). Then \eqref{eq:off_diagonal2@away} and the same balancing argument yield
	\[
		\norm[\big]{T^\lambda_J(T^\lambda_{J^\prime})^*}_{\op}
		\lesssim \lambda^{-\frac{1}{\delta_0(\lambda)+1}} 2^{-\frac{1}{2}|j_1-j_1^\prime|}.
	\]
Since \(|j_2-j_2^\prime| \le 2\), the bound \eqref{eq:CS_claim2@away} follows.

Applying the Cotlar--Stein lemma to \(\set[\big]{T^\lambda_J}_{J\in\Z^2 \cap \fm_{\textrm{bot}}^\vee}\), we obtain \(\norm[\big]{T^\lambda_\mathrm{ bot}}_{\op} \lesssim \lambda^{-\frac{1}{2(\delta_0(\lambda)+1)}}\).
The same argument with \(\fm_{\textrm{top}}\) in place of \(\fm_{\textrm{bot}}\) gives the same bound for \(T^\lambda_\mathrm{ top}\), and hence for \(T^\lambda_{0,s}\) uniformly in \(s\). Summing over \(0\le s \le K\) gives
	\begin{equation} \label{eq:UB_T0@away}
		\norm{T^\lambda_0}_{\op}
		\lesssim \lambda^{-\frac{1}{2(\delta_0(\lambda)+1)}}.
	\end{equation}


\subsection{Near the zero curves} \label{subsec:UB_near}

Fix \(1 \le \ell \le K\). By \eqref{eq:supp_fts}, the kernel of \(T^\lambda_\ell\) is supported where \(|y-c_\ell x^\gamma| \le \epsilon_0 x^\gamma\). We establish
	\[
		\norm{T^\lambda_\ell}_{\op} \lesssim \lambda^{-\frac{1}{2(\delta_\ell(\lambda)+1)}},
	\]
where \(\delta_\ell(\lambda)\) is defined in \eqref{eq:delta_l}.
We first treat the case \(\gamma \neq 0\); we defer the horizontal case \(\gamma = 0\) to the end of this subsection.


\subsubsection{Non-horizontal curves: decomposition by admissible boxes}

We first decompose \(T^\lambda_\ell\) dyadically in \(x\):
	\[
		T^\lambda_\ell = \sum_{j_1\in \Z} T^\lambda_{\ell,j_1}, \quad \text{where }
		T^\lambda_{\ell,j_1} f(x)
		= \int \ee^{\ii\lambda P(x,y)} \psi_\ell (x,y) \chi\paren[\Big]{\frac{x}{2^{-j_1}}} f(y) \dy.
	\]
The kernel of  \(T^\lambda_{\ell,j_1}\) is supported where \(x\approx 2^{-j_1}\) and \(y \approx 2^{-\gamma j_1}\).
Since \(\gamma\neq0\), the \(x\)- and \(y\)-projections of these supports have uniformly bounded overlap as \(j_1\) varies. Let \(\widetilde \chi\) be a slightly enlarged cutoff such that \(\widetilde\chi(y/2^{-\gamma j_1})\equiv 1\) on the \(y\)-support of \(T^\lambda_{\ell,j_1}\).
Setting \(f_{j_1}(y) \defeq \widetilde\chi\paren[\big]{\frac{y}{2^{-\gamma j_1}}} f(y)\), one has \(T^\lambda_{\ell,j_1}f=T^\lambda_{\ell,j_1}f_{j_1}\), and therefore
	\[
		\norm[\bigg]{\sum_{j_1\in \Z} T^\lambda_{\ell,j_1}f}_2^2
		\lesssim
		\sum_{j_1\in \Z} \norm[\big]{T^\lambda_{\ell,j_1}f_{j_1}}_2^2 
		\le
		\sup_{j_1} \norm[\big]{T^\lambda_{\ell,j_1}}_{\op}^2
			\sum_{j_1\in \Z} \norm{f_{j_1}}_2^2  
		\lesssim
		\sup_{j_1} \norm[\big]{T^\lambda_{\ell,j_1}}_{\op}^2 \norm{f}_2^2.
	\]
Hence,
	\begin{equation} \label{eq:j1_reduction@near}
		\norm[\bigg]{\sum_{j_1\in \Z} T^\lambda_{\ell,j_1}}_{\op}
		\lesssim
		\sup_{j_1} \norm[\big]{T^\lambda_{\ell,j_1}}_{\op}.
	\end{equation}

Fix \(j_1\). To estimate \(T^\lambda_{\ell,j_1}\), we decompose it by the dyadic transverse distance to the curve and then by the \(x\)-partition from \Cref{prop:admissible_box}. Let \(I_{j_1}\) denote the dyadic \(x\)-interval on which \(\chi(x/2^{-j_1})\) is supported, and set \(h_\mu \defeq 2^{-\gamma j_1-\mu}\). By the definition of \(\chi\) and the choice of \(\epsilon_0\) in \eqref{eq:epsilon_separation}, one has 
	\[
		\sum_{\mu \ge 0} \chi\paren[\Big]{\frac{y-c_\ell x^\gamma}{h_\mu}} =1
	\]
almost everywhere on the support of \(\psi_\ell(x,y)\chi(x/2^{-j_1})\).

For each \(\mu\) contributing to the sum, the support conditions give \(h_\mu \le 2^{|\gamma|+1} \epsilon_0 2^{-\gamma j_1}\). Together with \eqref{eq:epsilon_box}, this verifies the hypothesis of \Cref{prop:admissible_box}. Hence, for each such \(\mu\), we use the partition \(\set[\big]{[x_\nu, x_{\nu+1}]}_\nu\) of \(I_{j_1}\) provided by \Cref{prop:admissible_box} with \(h=h_\mu\). We write 
	\[
		T^\lambda_{\ell,j_1}= \sum_{\mu\ge 0} \sum_{\nu} T^\lambda_{\mu, \nu},
	\]
where, suppressing the fixed indices \(\ell\) and \(j_1\) from the notation,
	\[
		T^\lambda_{\mu, \nu} f(x)
		= \int \ee^{\ii\lambda P(x,y)}
			\psi_\ell (x,y)\chi\paren[\Big]{\frac{x}{2^{-j_1}}}
			\chi\paren[\Big]{\frac{y-c_\ell x^\gamma}{h_\mu}}
			\mathds 1_{[x_\nu, x_{\nu+1}]}(x)
			f(y) \dy.
	\]
We split
	\[
		\chi\paren[\Big]{\frac{y-c_\ell x^\gamma}{h_\mu}} 
		= \chi_+\paren[\Big]{\frac{y-c_\ell x^\gamma}{h_\mu}} 
		+ \chi_-\paren[\Big]{\frac{y-c_\ell x^\gamma}{h_\mu}},
	\]
where \(\supp \chi_+ \subset [\frac12,2]\) and \(\supp \chi_- \subset [-2, -\frac12]\). The argument for \(\chi_-\) is identical, so we restrict attention to the \(\chi_+\)-part and continue to write \(T^\lambda_{\mu,\nu}\) for the resulting block. Its kernel is supported in the curved piece
	\[
		R_{\nu,h_\mu}
		= \set[\bigg]{ (x,y) : x \in [x_\nu, x_{\nu+1}], \ \frac{h_\mu}{2} \le y-c_\ell x^\gamma \le 2h_\mu}.
	\]
Thus \Cref{prop:admissible_box} gives the associated admissible box \(B_{\nu,h_\mu}\) with \(R_{\nu, h_\mu} \subset B_{\nu,h_\mu}\). Moreover, by \ref{item:A2} in \Cref{prop:admissible_box}, as \(\nu\) varies, the \(x\)- and \(y\)-projections of the boxes \(\set{B_{\nu, h_\mu}}_\nu\) have uniformly bounded overlap. Therefore,
	\[
		\norm[\bigg]{\sum_\nu T^\lambda_{\mu,\nu}}_{\op}
		\lesssim \sup_\nu \norm[\big]{T^\lambda_{\mu,\nu}}_{\op}.
	\]
Summing over \(\mu\), we obtain
	\begin{equation} \label{eq:reduce_mu_sum@near}
		\norm[\big]{T^\lambda_{\ell,j_1}}_{\op}
		\lesssim
		\sum_{\mu\ge0}
		\sup_\nu \norm[\big]{T^\lambda_{\mu,\nu}}_{\op}.
	\end{equation}


\subsubsection{Non-horizontal curves: block estimates}

To estimate a fixed block \(T^\lambda_{\mu,\nu}\), we introduce the adapted scale vector
	\[
		\tJ = (\tj_1, \tj_2) \defeq (j_1, \gamma j_1 + \mu).
	\]
By construction, one has
	\[
		x \approx 2^{-\tj_1} \quad \text{and} \quad Y \approx 2^{-\tj_2} \quad \text{on } B_{\nu, h_\mu}.
	\]
Since \(\bv_\ell = (1-\gamma,2)\) by \eqref{eq:NP_near}, 
the area estimate \ref{item:A1} in \Cref{prop:admissible_box} yields
	\[
		|B_{\nu,h_\mu}| \approx x^{1-\gamma}Y^2 \approx 2^{-\tJ\cdot \bv_\ell}.
	\]
Hence \Cref{lemma:size_est} yields the size estimate
	\[
		\norm[\big]{T^\lambda_{\mu,\nu}}_{\op}
		\lesssim 2^{-\tJ\cdot \bv_\ell/2}.
	\]
	
For the oscillatory estimate, recall from \eqref{eq:vertices_Pl} that \(\fm_\ell = (m_{\ell,1}, m_{\ell,2})\). Applying \Cref{lemma:admissible_box_hessian} to the box \(B_{\nu,h_\mu}\) gives
	\[
		|P_{xy}(x,y)| \approx  |x|^{m_{\ell,1}}|Y|^{m_{\ell,2}} \approx 2^{-\tJ\cdot \fm_\ell}
		\quad \text{for } (x,y) \in B_{\nu,h_\mu},
	\]
together with the corresponding \(y\)-derivative bounds at scale \(2^{-\tj_2}\).
The vertical side length of \(B_{\nu,h_\mu}\) is \(\approx h_\mu = 2^{-\tj_2}\), and the cutoff factors satisfy the corresponding \(y\)-derivative bounds. Since \(B_{\nu, h_\mu}\) is an axis-parallel box, every horizontal segment required in \Cref{lemma:osc_est} remains in \(B_{\nu,h_\mu}\). Hence, applying \Cref{lemma:osc_est} with \(\omega = 2^{-\tJ\cdot \fm_\ell}\) yields
	\[
		\norm[\big]{T^\lambda_{\mu,\nu}}_{\op}
		\lesssim 
		\lambda^{-1/2}2^{\tJ\cdot\fm_\ell/2}.
	\]
Combining these two estimates, we obtain
	\begin{equation}\label{eq:block_mu_tildeJ@near}
		\norm[\big]{T^\lambda_{\mu,\nu}}_{\op}
		\lesssim
		\min\set[\Big]{2^{-\tJ\cdot\bv_\ell/2}, \lambda^{-1/2}2^{\tJ\cdot\fm_\ell/2}}.
	\end{equation}


\subsubsection{Non-horizontal curves: summation in the transverse scale}

To sum the block estimates over \(\mu\), we distinguish two cases. Assume first that \(\delta_\ell(\lambda)<+\infty\). By \Cref{lemma:NP_geometry}, one has \(\P_\ell = \fm_\ell + \cone\set{(-\gamma,1)}\). Then by definition,
	\[
		\fm_\ell^\vee 
		= \set{\fq \in \R^2 : \fq \cdot \fm_\ell = \min_{\fp \in \P_\ell} \fq \cdot \fp}
		= \set{\fq \in \R^2: \fq\cdot (-\gamma,1) \ge 0}.
	\]
Since \(\tJ \cdot (-\gamma,1) = \mu \ge 0\), one has \(\tJ \in \fm_\ell^\vee\). As \(\delta_\ell(\lambda)\bv_\ell \in \P_\ell\), it follows that \(\tJ \cdot \fm_\ell \le \delta_\ell(\lambda)\tJ \cdot \bv_\ell\). Substituting this into \eqref{eq:block_mu_tildeJ@near} gives
	\[
		\norm[\big]{T^\lambda_{\mu,\nu}}_{\op}
		\lesssim \min 
		\set[\Big]{
			2^{-\tJ\cdot\bv_\ell/2}, \lambda^{-1/2}2^{\delta_\ell(\lambda)\tJ\cdot \bv_\ell/2}
		}.
	\]
By \Cref{cor:ray_geometry}, one has \(\delta_\ell(\lambda)>0\). 
Since \(\tJ\cdot \bv_\ell = (1+\gamma)j_1 + 2\mu\), splitting the sum over \(\mu\) at \(2^{-\tJ \cdot \bv_\ell} \approx \lambda^{-1/(\delta_\ell(\lambda)+1)}\) yields
	\[
		\sum_{\mu \ge 0}
		\sup_\nu \norm[\big]{T^\lambda_{\mu,\nu}}_{\op}
		\lesssim \lambda^{-\frac{1}{2(\delta_\ell(\lambda)+1)}}.
	\]
Combining this estimate with \labelcref{eq:j1_reduction@near,eq:reduce_mu_sum@near}, we obtain
	\[
		\norm[\big]{T^\lambda_\ell}_{\op} \lesssim \lambda^{-\frac{1}{2(\delta_\ell(\lambda)+1)}}.
	\]

We next consider the case \(\delta_\ell(\lambda) = +\infty\). By \Cref{cor:ray_geometry}, this can occur only when \(\gamma = -1\). Then we have \(\bv_\ell = (2,2)\) and \(\tJ=(j_1, -j_1+\mu)\), so that \(\tJ\cdot \bv_\ell =2\mu\).
Substituting the size estimate \(\norm{T^\lambda_{\mu,\nu}}_{\op} \lesssim 2^{-\mu}\) from the block estimate \eqref{eq:block_mu_tildeJ@near} into \eqref{eq:reduce_mu_sum@near}, we obtain
	\[
		\norm[\big]{T^\lambda_{\ell,j_1}}_{\op} \lesssim 1 \text{ uniformly in } j_1.
	\] 
The reduction \eqref{eq:j1_reduction@near} then gives
	\[
		\norm[\big]{T^\lambda_\ell}_{\op}\lesssim 1 = \lambda^{-\frac{1}{2(\delta_\ell(\lambda)+1)}}.
	\]
	
Consequently, for every \(\gamma \neq 0\),
	\begin{equation} \label{eq:Tl_bound_gamma_nonzero@near}
		\norm[\big]{T^\lambda_\ell}_{\op} \lesssim \lambda^{-\frac{1}{2(\delta_\ell(\lambda)+1)}}.
	\end{equation}


\subsubsection{Horizontal curves}

Finally, assume \(\gamma=0\), and write \(Y=y-c_\ell\). Since \(Y\) differs from \(y\) only by a translation, ordinary dyadic blocks in the adapted coordinates \((x,Y)\) are also axis-parallel boxes in the original coordinates \((x,y)\). We decompose
	\[
		T^\lambda_\ell = \sum_{j_1 \in \Z} \sum_{\mu \ge 0} T^\lambda_{j_1, \mu},
		\quad \text{where }
		T^\lambda_{j_1,\mu} f(x)
		\defeq \int \ee^{\ii\lambda P(x,y)}\psi\paren[\Big]{\frac{y-c_\ell}{\epsilon_0}}
			\chi\paren[\Big]{\frac{x}{2^{-j_1}}}\chi\paren[\Big]{\frac{y-c_\ell}{2^{-\mu}}}f(y)\dy.
	\]
As before, we restrict attention to the part supported where \(Y>0\) and continue to write \(T^\lambda_{j_1,\mu}\) for the resulting block; the other part supported where \(Y<0\) is treated identically.

Set \(\tJ=(\tj_1,\tj_2)\defeq(j_1, \mu)\). The kernel is supported in an axis-parallel box in the original coordinates on which
	\[
		x \approx 2^{-\tj_1} \quad \text{and} \quad Y=y-c_\ell \approx 2^{-\tj_2}.
	\]
Since \(\bv_\ell = (1,1)\) by \eqref{eq:NP_near}, this box has area \(\approx  2^{-\tJ\cdot \bv_\ell}\). Moreover, on the support of the block, the choice of \(\epsilon_0\) gives \(|y| \approx 1\) and \(|y-c_n| \approx 1\) for every \(n \neq \ell\). Combining these bounds with the factorization
	\[
		P_{xy}(x,y)=ux^My^N \prod_{n=1}^L (y-c_n)^{k_n}
	\]
and \eqref{eq:vertices_Pl}, we obtain
	\[
		|P_{xy}(x,y)| \approx |x|^{m_{\ell,1}}|Y|^{m_{\ell,2}}
		\approx 2^{-\tJ\cdot \fm_\ell}.
	\]
Since \(\partial_y = \partial_Y\), differentiating the same factorization gives, for \(b=1,2\),
	\[
		|\partial_y^b P_{xy}(x,y)| \lesssim |x|^{m_{\ell,1}}|Y|^{m_{\ell,2}-b} \approx 2^{-\tJ\cdot\fm_\ell}(2^{-\tj_2})^{-b}.
	\]
The vertical side length of the box is \(\approx 2^{-\tj_2}\), and the cutoff factors satisfy the corresponding \(y\)-derivative bounds. Since the box is axis-parallel, every horizontal segment required in \Cref{lemma:osc_est} remains in the box. Applying \Cref{lemma:size_est,lemma:osc_est} then yields
	\begin{equation} \label{eq:block_j1_mu_tildeJ}
		\norm[\big]{T^\lambda_{j_1,\mu}}_{\op}
		\lesssim
		\min\set[\Big]{
			2^{-\tJ\cdot \bv_\ell/2}, \lambda^{-1/2}2^{\tJ\cdot\fm_\ell/2}
		}.
	\end{equation}
	
For \(\gamma =0\), \Cref{lemma:NP_geometry} gives \(\P_\ell = \fm_\ell + \cone\set{(0,1)}\). Then by definition,
	\[
		\fm_\ell^\vee
		= \set{\fq \in \R^2 : \fq \cdot \fm_\ell = \min_{\fp \in \P_\ell} \fq \cdot \fp}
		= \set{\fq \in \R^2: \fq \cdot (0,1) \ge 0}.
	\]
 Since \(\tJ\cdot(0,1)=\mu \ge0\), one has \(\tJ \in \fm_\ell^\vee\). As \(\delta_\ell(\lambda)\bv_\ell \in \P_\ell\), it follows that \(\tJ\cdot \fm_\ell \le \delta_\ell(\lambda) \tJ\cdot \bv_\ell\). Substituting this into \eqref{eq:block_j1_mu_tildeJ} gives
	\[
		\norm[\big]{T^\lambda_{j_1,\mu}}_{\op}
		\lesssim \min
		\set[\Big]{2^{-\tJ\cdot \bv_\ell/2}, \lambda^{-1/2}2^{\delta_\ell(\lambda)\tJ\cdot\bv_\ell/2}}.
	\]
We sum these blocks along the diagonals \(d=\tJ\cdot\bv_\ell=j_1+\mu\). For a fixed \(d\), the kernels of \(\set{T^\lambda_{j_1, d-j_1}}_{j_1}\) are supported where \(x \approx 2^{-j_1}\) and \(Y \approx 2^{-d+j_1}\). Their \(x\)- and \(Y\)-projections thus have uniformly bounded overlap. Since \(Y=y-c_\ell\), the same holds for their \(x\)- and \(y\)-projections in the original coordinates. Therefore,
	\[
		\norm[\bigg]{\sum_{\substack{j_1+\mu=d \\ \mu \ge 0}} T^\lambda_{j_1,\mu}}_{\op}
		\lesssim \sup_{\substack{j_1+\mu=d \\ \mu \ge 0}} \norm[\big]{T^\lambda_{j_1,\mu}}_{\op}
		\lesssim \min\set[\Big]{2^{-d/2}, \lambda^{-1/2}2^{\delta_\ell(\lambda)d/2}}.
	\]
Since \(\delta_\ell(\lambda) >0\) by \Cref{cor:ray_geometry}, we may sum over \(d \in \Z\), splitting at \(2^{-d} \approx \lambda^{-1/(\delta_\ell(\lambda)+1)}\), to obtain
	\begin{equation} \label{eq:Tl_bound_gamma_zero@near}
		\norm{T^\lambda_\ell}_{\op} 
		\lesssim \lambda^{-\frac{1}{2(\delta_\ell(\lambda)+1)}}.
	\end{equation}

Combining \labelcref{eq:Tl_bound_gamma_nonzero@near,eq:Tl_bound_gamma_zero@near}, we conclude that
	 \begin{equation} \label{eq:UB_Tl@near}
	 	\norm[\big]{T^\lambda_\ell}_{\op}
		\lesssim
		\lambda^{-\frac{1}{2(\delta_\ell(\lambda)+1)}}
		\quad \text{for every } 1 \le \ell \le K.
	 \end{equation}
Finally, summing the estimates \labelcref{eq:UB_T0@away,eq:UB_Tl@near} and using the definitions of \(\df(P)\) and \(\db(P)\) in \eqref{eq:global_newton_distances}, we obtain
	\[
		\norm[\big]{T^\lambda}_{\op}
		\lesssim
		\begin{cases}
			\lambda^{-\frac{1}{2(\df(P)+1)}}, & \text{if } \lambda \ge 1, \\[1ex]
			\lambda^{-\frac{1}{2(\db(P)+1)}}, & \text{if } 0< \lambda <1.
		\end{cases}
	\]
This proves the upper bound in \Cref{thm:main}\,(A).


\section{Lower bounds for balanced Newton pairs} \label{sec:lower_bounds}

We now prove the matching lower bound in \Cref{thm:main}\,(A), assuming that the family of Newton pairs \(\fN(P)\), defined in \eqref{eq:NF}, is balanced. For each balanced Newton pair for which the relevant Newton distance is finite, we construct a box on which the phase oscillation is negligible and whose area realizes that Newton distance.

We begin with a standard Knapp-type test on rectangular boxes. We will also use this lemma in \Cref{sec:unboundedness}.

\Lemma{Knapp-type box test}{\label{lemma:box_testing}
Let
	\[
		Q = \set{(x,y): |x-x_0| \le \Delta x, \ |y-y_0| \le \Delta y} \subset \R_{>0}^2.
	\]
If \(\lambda |Q| \sup_Q |P_{xy}| \le 1\), then
	\[
		\sup_{0\neq f\in\cS(\R)} \frac{\norm{T^\lambda f}_2}{\norm{f}_2} \ge \frac{|Q|^{\frac12}}{4}.
	\]
}
\pf{
By the fundamental theorem of calculus, for \((x,y) \in Q\), we have
	\[
		\begin{aligned}
			\lambda\abs{P(x,y)-P(x_0,y)-P(x,y_0)+P(x_0,y_0)} 
			&= \lambda\abs[\bigg]{\int_{x_0}^x \int_{y_0}^y P_{xy}(u,v) \d{v} \d{u}} \\
			&\le \lambda \Delta x \Delta y \sup_Q |P_{xy}| \le \frac{1}{4},
		\end{aligned}
	\]
and hence \(\Re\paren[\big]{\ee^{\ii\lambda(P(x,y)-P(x_0,y)-P(x,y_0)+P(x_0,y_0))}} \ge \frac{1}{2}\) on \(Q\).
We define the test functions
	\[
		f(y) = \ee^{-\ii\lambda P(x_0, y)}\psi\paren[\Big]{\frac{y-y_0}{\Delta y}}, \qquad
		g(x) = \ee^{\ii\lambda(P(x,y_0)-P(x_0,y_0))}\psi\paren[\Big]{\frac{x-x_0}{\Delta x}},
	\]
so that \(\supp(g \otimes f) \subset Q\). Since \(\psi^2 \le \psi\),
	\[
		\abs{\inn{T^\lambda f, g}} 
		\ge \Re \inn{T^\lambda f, g} 
		\ge \frac{1}{2} 
			\iint \psi\paren[\Big]{\frac{y-y_0}{\Delta y}}
			\psi\paren[\Big]{\frac{x-x_0}{\Delta x}}
			\d{y}\d{x}
		\ge \frac12 \norm{f}_2^2\norm{g}_2^2.
	\]
Moreover, one has \(\norm{f}_2\norm{g}_2 \ge |Q|^{1/2}/2\). Therefore
	\[
		\sup_{0\neq h\in\cS(\R)} \frac{\norm{T^\lambda h}_2}{\norm{h}_2}
		\ge \frac{\abs{\inn{T^\lambda f, g}}}{\norm{f}_2\norm{g}_2}
		\ge \frac{|Q|^{\frac{1}{2}}}{4}. \qedhere
	\]
}

We now apply \Cref{lemma:box_testing} to boxes adapted to a fixed balanced Newton pair. Fix \(\lambda>0\) and a constant \(0<\epsilon\ll 1\), independent of \(\lambda\).

Recall the definition of \(\delta_\ell(\lambda)\) from \eqref{eq:delta_l}. Since the balanced family \(\fN(P)\) is finite, the extrema in \eqref{eq:global_newton_distances} are attained. We therefore choose an index \(0 \le \ell \le K\) such that
	\[
		\delta_\ell(\lambda)=
		\begin{cases}
			\df(P), & \text{if } \lambda \ge 1, \\[0.5ex]
			\db(P), & \text{if } 0< \lambda <1.
		\end{cases}
	\]
Here \(\delta_\ell(\lambda)\) is finite: every forward Newton distance is finite, while \(\db(P) \le \db(\P_0, \bv_0) <+\infty\). We shall prove
	\begin{equation} \label{eq:LB}
		\norm{T^\lambda}_{\op} \gtrsim \lambda^{-\frac{1}{2(\delta_\ell(\lambda)+1)}}.
	\end{equation}
By \Cref{lemma:box_testing}, it is enough to construct boxes \(Q_\alpha\) so that
	\[
		\lambda |Q_\alpha|\,\sup_{Q_\alpha} |P_{xy}| \le 1
		\qquad \text{and} \qquad
		|Q_\alpha| \approx \lambda^{-\frac{1}{\delta_\ell(\lambda)+1}}.
	\]

For the fixed index \(0 \le \ell \le K\), we use the exponent set \(\Lambda_\ell\) from \eqref{eq:Lambda_sets} and the vector \(\bv_\ell\) from \labelcref{eq:NP_away,eq:NP_near}. We construct \(\fq \in \R^2\) and \(\alpha \in \R\) satisfying
	\begin{equation} \label{eq:choice_alpha@LB}
		2^{-\alpha(\fq\cdot\bv_\ell)} = \lambda^{-\frac{1}{\delta_\ell(\lambda)+1}}
	\end{equation}
and
	\begin{equation} \label{eq:supporting_condition@LB}
		2^{-\alpha(\fq\cdot\fm)} \le 2^{-\alpha(\fq\cdot\bv_\ell)\delta_\ell(\lambda)}
		\quad \text{for all } \fm \in \Lambda_\ell.
	\end{equation}
By \Cref{cor:ray_geometry}, the construction of \(\fq\) and \(\alpha\) splits into the cases \(\gamma \neq -1\) and \(\gamma = -1\).

\paragraph{The case \(\gamma \neq -1\).}
We set \(\fq=(1,\gamma)\). By \Cref{lemma:NP_geometry}, \(\P_\ell\) is then contained in the affine line \(\plane{\fq}{r}=\set[\big]{\by\in \R^2: \fq \cdot \by = r}\) for some \(r \in \R\). Since \((\P_\ell, \bv_\ell)\) is balanced, one has \(\delta_\ell(\lambda)\bv_\ell \in \P_\ell\). Hence, for every \(\fm \in \Lambda_\ell \subset \P_\ell\), we have 
	\[
		\fq\cdot \fm = (\fq\cdot\bv_\ell)\delta_\ell(\lambda).
	\]
Since \(\fq\cdot \bv_\ell=1+\gamma \neq0\), there exists \(\alpha \in \R\) satisfying \eqref{eq:choice_alpha@LB}. Then \eqref{eq:supporting_condition@LB} holds with equality.

\paragraph{The case \(\gamma = -1\).}
By \Cref{cor:ray_geometry}, the balanced pair \((\P_\ell, \bv_\ell)\) satisfies 
	\[
		\P_\ell = [\df(\P_\ell, \bv_\ell), \db(\P_\ell, \bv_\ell)]\bv_\ell.
	\]
We set \(\fq = \bv_\ell\). The line through the forward endpoint with inward normal \(\fq\) supports \(\P_\ell\); if \(\db(\P_\ell,\bv_\ell)<+\infty\), the line through the backward endpoint with inward normal \(-\fq\) also supports \(\P_\ell\). Equivalently, for every \(\fm \in \Lambda_\ell\),
	\[
		\df(\P_\ell,\bv_\ell)\,(\fq \cdot \bv_\ell) \le \fq \cdot \fm \le \db(\P_\ell,\bv_\ell)\,(\fq\cdot\bv_\ell).
	\]
Since \(\fq\cdot \bv_\ell>0\), there exists \(\alpha\in\R\) satisfying \eqref{eq:choice_alpha@LB}. For this choice of \(\alpha\), we obtain \eqref{eq:supporting_condition@LB} as follows:
	\begin{itemize}
		\item
			if \(\lambda \ge 1\), then \(\alpha \ge 0\) and \(\delta_\ell(\lambda)=\df(\P_\ell,\bv_\ell)\), so the first inequality yields \eqref{eq:supporting_condition@LB};
		\item
			if \(0<\lambda<1\), then \(\alpha <0\) and \(\delta_\ell(\lambda)=\db(\P_\ell,\bv_\ell)\), so the second inequality yields \eqref{eq:supporting_condition@LB}.
	\end{itemize}
\medskip
 
In either case, write \(\fq=(q_1,q_2)\) and construct the box \(Q_\alpha\) as follows. If \(\ell = 0\), set
	\[
		x_\alpha = 2^{-\alpha q_1}, \qquad
		y_\alpha = 2^{-\alpha q_2}, \qquad 
		\Delta x_\alpha = \epsilon x_\alpha, \qquad 
		\Delta y_\alpha = \epsilon y_\alpha.
	\]
If \(1 \le \ell \le K\), set
	\[
		\begin{alignedat}{2}
			&x_\alpha = 2^{-\alpha q_1}, \qquad
			Y_\alpha = 2^{-\alpha q_2}, \qquad
			&&y_\alpha = c_\ell x_\alpha^\gamma+Y_\alpha, \\[0.5ex]
			&\Delta x_\alpha =
				\begin{cases}
					\epsilon x_\alpha^{1-\gamma}Y_\alpha, & \text{if } \gamma \neq 0, \\[0.5ex]
					\epsilon x_\alpha, & \text{if } \gamma = 0,
				\end{cases}
			\qquad
			&&\Delta y_\alpha = \epsilon Y_\alpha.
		\end{alignedat}
	\]
In both cases, we define 
	\[
		Q_\alpha \defeq \set[\big]{(x,y): |x-x_\alpha| \le \Delta x_\alpha, \ |y-y_\alpha| \le \Delta y_\alpha}.
	\]

Before applying \Cref{lemma:box_hessian}, we verify the additional hypothesis \(x_\alpha^{-\gamma}Y_\alpha \le 1\) of \ref{item:T2} in \Cref{lemma:box_hessian}, which is relevant only when \(1 \le \ell \le K\) and \(\gamma \neq 0\).
If \(1 \le \ell \le K\) and \(\gamma \neq 0\),
 	\[
 		x_\alpha^{-\gamma}Y_\alpha = 2^{-\alpha \fq \cdot (-\gamma,1)}=
 		\begin{cases}
 			2^{-\alpha(1,\gamma)\cdot(-\gamma,1)} =1, & \text{if } \gamma \neq -1, \\[0.5ex]
 			2^{-\alpha(2,2)\cdot(1,1)} = 2^{-4\alpha}, & \text{if } \gamma = -1.
 		\end{cases}
 	\]
 In the latter case \(\gamma = -1\), \Cref{cor:ray_geometry} yields \(\db(\P_\ell,\bv_\ell)=+\infty\). Since \(\delta_\ell(\lambda)\) is finite, this forces \(\lambda \ge 1\), and hence \(\alpha \ge 0\) by \eqref{eq:choice_alpha@LB}. Therefore \(x_\alpha^{-\gamma}Y_\alpha \le 1\) in all cases, as required.
 
 We now apply \Cref{lemma:box_hessian}, using \ref{item:T1} when \(\ell = 0\) and \ref{item:T2} when \(1 \le \ell \le K\). By \eqref{eq:supporting_condition@LB},
 	\[
 		\sup_{Q_\alpha} |P_{xy}| 
 		\lesssim \max_{\fm \in \Lambda_\ell} 2^{-\alpha(\fq \cdot \fm)}
 		\le 2^{-\alpha(\fq \cdot \bv_\ell)\delta_\ell(\lambda)}.
 	\]
 Moreover, the definition of \(\bv_\ell\) in \labelcref{eq:NP_away,eq:NP_near} ensures \(|Q_\alpha| \approx \epsilon^2 2^{-\alpha(\fq\cdot\bv_\ell)}\) in all cases. Therefore, by \eqref{eq:choice_alpha@LB},
 	\[
		\lambda|Q_\alpha| \sup_{Q_\alpha} |P_{xy}|
		\lesssim \epsilon^2 \lambda 2^{-\alpha(\fq\cdot\bv_\ell)(\delta_\ell(\lambda)+1)} 
		\lesssim \epsilon^2.
	\]
For sufficiently small \(\epsilon>0\), \Cref{lemma:box_testing} and the choice of \(\alpha\) in \eqref{eq:choice_alpha@LB} yield
	\[
		\norm{T^\lambda}_{\op} 
		\gtrsim |Q_\alpha|^{\frac 12}
		\approx 2^{-\alpha(\fq\cdot\bv_\ell)/2}
		= \lambda^{-\frac{1}{2(\delta_\ell(\lambda)+1)}}.
	\]

Consequently, by the choice of \(\ell\),
	\[
		\norm{T^\lambda}_{\op} \gtrsim
		\begin{cases}
			\lambda^{-\frac{1}{2(\df(P)+1)}}, & \text{if } \lambda \ge 1, \\[1ex]
			\lambda^{-\frac{1}{2(\db(P)+1)}}, & \text{if } 0< \lambda <1.
		\end{cases}
	\]
This proves the lower bound in \Cref{thm:main}\,(A).
Together with the upper bounds proved in \Cref{sec:upper_bounds}, this completes the proof of \Cref{thm:main}\,(A).


\section{Unboundedness for unbalanced Newton pairs} \label{sec:unboundedness}

We prove \Cref{thm:main}\,(B). Fix \(\lambda>0\). We also fix a constant \(0<\epsilon\ll1\) not depending on \(\lambda\), chosen sufficiently small for all boxes below.

We first introduce an elementary separation lemma needed for the box construction.

\Lemma{Supporting-line separation for unbalanced pairs}{\label{lemma:unbalancing_criterion}
Let \(\Lambda \subset \R \times \R_{\ge0}\) and \(B \subset \Q^2\) be nonempty finite sets. For \(\P=\conv\paren[\big]{\Lambda+\cone B}\), assume that \(\P \subset \R \times \R_{\ge 0}\). Let \(\bv \in \R \times \R_{>0}\), and suppose
	\[
		\P \cap \cone(\bv) = \emptyset.
	\]
Then there exist \(\fq \in \R^2 \setminus \set{\bm 0}\) and \(r>0\) such that
	\begin{enumerate}[label=\rm (\arabic*)]
		\item \(\fq \cdot \bv=0\),
		\item \(\fq \cdot \fm \ge r\) for all \(\fm \in \Lambda\),
		\item \(\fq \cdot \fb \ge 0\) for all \(\fb \in B\).	
	\end{enumerate}
Consequently, \(\fq\cdot \fp \ge r\) for all \(\fp \in \P\).
}
\pf{
The assumptions imply \(\P \cap \R\bv = \emptyset\). Since \(\P\) is convex, it is contained in one of the two open half-planes bounded by \(\R\bv\). Hence we may choose a normal vector \(\fq\) to \(\R\bv\) such that \(\fq \cdot \fp>0\) for every \(\fp \in \P\); in particular \(\fq\cdot\bv=0\), which proves (1).
Since \(\Lambda \subset \P\) and \(\Lambda\) is finite, we may set \(r= \min_{\fm \in \Lambda} \fq \cdot \fm >0\), which proves (2).
Finally, since \(\fp+t\fb \in \P\) for every \(\fp \in \P\), \(\fb \in B\) and \(t\ge0\), the functional \(\fq\) would be unbounded below on \(\P\) unless \(\fq \cdot \fb \ge 0\), which proves (3). Since \(\P=\conv(\Lambda + \cone B)\), (2) and (3) yield \(\fq \cdot \fp \ge r\) for every \(\fp \in \P\), as claimed.
}

Assume now that the family \(\fN(P)\), defined in \eqref{eq:NF}, is unbalanced. We distinguish two subcases:
	\begin{enumerate}[label=(\roman*)]
		\item 
			the pair \((\P_0, \bv_0)\) is unbalanced;
		\item
			the pair \((\P_0, \bv_0)\) is balanced, but \((\P_\ell, \bv_\ell)\) is unbalanced for some \(1 \le \ell \le K\).
	\end{enumerate}
In particular, \Cref{cor:ray_geometry} gives \(\gamma \neq -1\) in case (ii).

In case (i), set \(\ell=0\). In case (ii), fix \(1\le\ell\le K\) such that \((\P_\ell, \bv_\ell)\) is unbalanced. Recall from \labelcref{eq:NP_away,eq:NP_near} that 
	\[
		\begin{alignedat}{2}
			&\P_0 &&=\conv\paren[\big]{\Lambda_0 + \cone\set{\bm 0}}, \\[0.5ex]
			&\P_\ell &&= \conv\paren[\big]{\Lambda_\ell + \cone\set{(-\gamma,1)}}
			\quad \text{for } 1 \le \ell \le K.
		\end{alignedat}
	\]
Hence, applying \Cref{lemma:unbalancing_criterion} to the unbalanced pair gives \(\fq \in \R^2\setminus\set{\bm 0}\) and \(r>0\) such that
	\[
		\fq\cdot\bv_\ell = 0
		\qquad \text{and} \qquad
		\fq\cdot\fm \ge r \quad \text{for all } \fm \in \Lambda_\ell.
	\]
In case (ii), \Cref{lemma:unbalancing_criterion} also gives \(\fq \cdot (-\gamma,1) \ge 0\). Since \(\gamma \neq -1\) in case (ii), the vector \((-\gamma,1)\) is not parallel to \(\bv_\ell\). Hence the orthogonality  \(\fq \cdot \bv_\ell = 0\) forces \(\fq \cdot (-\gamma,1) \neq 0\), and thus \(\fq\cdot(-\gamma,1)>0\).

To apply \Cref{lemma:box_testing} and ensure
	\[
		|Q_\alpha| \toinf \quad \text{while} \quad \lambda|Q_\alpha|\sup_{Q_\alpha} |P_{xy}| \to 0 \quad \text{as } \alpha \toinf,
	\]
we perturb \(\fq\) slightly in the negative \(\bv_\ell\)-direction. For any \(\rho>0\), the perturbed vector \(\widetilde \fq \defeq \fq-\rho \bv_\ell\) satisfies
	\begin{equation} \label{eq:tildeq_common1@unbal}
		\widetilde \fq \cdot \bv_\ell = -\rho|\bv_\ell|^2 <0.
	\end{equation}
Since \(\Lambda_\ell\) is finite, we may choose \(\rho\) small enough that
	\begin{align}
	&\widetilde \fq \cdot \fm \ge \frac r2 \quad \text{for all } \fm \in \Lambda_\ell, \label{eq:tildeq_common2@unbal} \\
	&\rho|\bv_\ell|^2 < \frac r2.
		\label{eq:tildeq_common3@unbal}
	\end{align}
In case (ii), we also preserve \(\widetilde \fq \cdot (-\gamma,1)>0\), decreasing \(\rho\) if necessary.

For \(\alpha \gg 1\), we use the same box construction as in \Cref{sec:lower_bounds}, with \(\widetilde \fq\) replacing \(\fq\). In case (ii), if \(\gamma \neq 0\), then \(\widetilde \fq \cdot (-\gamma,1)>0\). This yields 
	\(
		x_\alpha^{-\gamma}Y_\alpha = 2^{-\alpha \widetilde \fq \cdot (-\gamma,1)} \le 1
	\)
for all sufficiently large \(\alpha\), so the additional hypothesis in \ref{item:T2} of \Cref{lemma:box_hessian} is satisfied.

We now estimate the area of \(Q_\alpha\) and the Hessian on \(Q_\alpha\). Applying \Cref{lemma:box_hessian}, using \ref{item:T1} in case (i) and \ref{item:T2} in case (ii), and then using \eqref{eq:tildeq_common2@unbal}, we obtain the Hessian bound
	\[
		\sup_{Q_\alpha} |P_{xy}|
		\lesssim \max_{\fm \in \Lambda_\ell} 2^{-\alpha \widetilde \fq \cdot \fm} \le 2^{-\alpha r/2}.
	\]
Meanwhile, the box construction and \eqref{eq:tildeq_common1@unbal} give
	\[
		|Q_\alpha| 
		\approx \epsilon^2 2^{-\alpha \widetilde \fq \cdot \bv_\ell}
		= \epsilon^2 2^{\alpha \rho |\bv_\ell|^2}.
	\]
Combining the bounds for \(\sup_{Q_\alpha} |P_{xy}|\) and \(|Q_\alpha|\) with \eqref{eq:tildeq_common3@unbal}, we obtain
	\[
		\lambda |Q_\alpha| \sup_{Q_\alpha} |P_{xy}|
		\lesssim \lambda \epsilon^2 2^{\alpha(\rho|\bv_\ell|^2-r/2)} \to 0
		\quad \text{as } \alpha \toinf.
	\]
Thus, for all sufficiently large \(\alpha\), \(\lambda |Q_\alpha| \sup_{Q_\alpha}|P_{xy}| \le 1\). Hence \Cref{lemma:box_testing} gives
	\[
		\sup_{0\neq f\in\cS(\R)} \frac{\norm{T^\lambda f}_2}{\norm{f}_2} 
		\gtrsim |Q_\alpha|^{1/2} \approx 2^{\alpha \rho |\bv_\ell|^2/2}.
	\]
Letting \(\alpha\toinf\) completes the proof of \Cref{thm:main}\,(B).


\section{Concluding remarks} \label{sec:conclusion}

The central point of this paper is that the pairing of each Newton polygon \(\P_\ell\) with the vector \(\bv_\ell\) recording the box-area scale determines the global \(L^2\)-estimates. For quasi-homogeneous phases, a single shear along each real zero curve \(y=c_\ell x^\gamma\) puts the Hessian into normal-crossing form, and all the polygons \(\P_\ell\) lie on a common line. Consequently, the sharp exponents are computed purely algebraically from the factorization of \(P_{xy}\), even when scaling does not reduce the global problem to the local theory.

For general polynomial phases, distinct Puiseux expansions of the real zero set of \(P_{xy}\) near infinity may share several initial terms, so a single shear need not put the Hessian into normal-crossing form. In a forthcoming paper, we organize the successive shears by finite trees of common partial sums and show that the Newton pairs attached to the nodes of these trees, together with the local Newton distances at finite points, characterize global \(L^2\)-boundedness and determine the sharp \(L^2\)-estimates for all real polynomial phases with \(P_{xy}\not\equiv0\). Since the polygons need not lie on a common line, supporting lines and dual faces play a central role there.

Further natural questions include the sharp \(L^p\)-\(L^q\) improving range of the global Radon transform, a global counterpart of the local estimates of Phong and Stein, and a cutoff-free \(L^2\)-classification for oscillatory integral operators in higher dimensions.


\section*{Acknowledgments}

Both authors were supported by the National Research Foundation of Korea (NRF) grant funded by the Korea government (MSIT) under Grant No.~RS-2026-25482410.


\appendix
\section{Factorization of quasi-homogeneous polynomials} \label{appendix:form_quasi}

In this appendix, we give a self-contained proof of the factorization \eqref{eq:quasi_form}.

\Lemma{}{\label{lemma:form_quasi}
Let \(P\in \R[x,y]\) be a quasi-homogeneous polynomial and assume \(P_{xy} \not \equiv 0\). Then, up to exchanging \(x\) and \(y\), \(P_{xy}\) can be factored for \(x>0\) as
	\[
		P_{xy}(x,y) = u x^M y^N \prod_{\ell=1}^{L}(y-c_\ell x^\gamma)^{k_\ell},
	\]
where \(u \in \R\setminus\set{0}\), \(M,N \in \Z_{\ge 0}\), \(k_\ell \in \N\), \(\gamma\in \Q\), and
\(c_1,\dots,c_L\in\C\setminus\set{0}\) are distinct.
}

\pf{
By quasi-homogeneity, the exponents \((m,n) \in \Lambda(P_{xy})\) lie on a line \(w_1 m+ w_2 n =d\) for some \(w_1, w_2, d \in \Z\) with \((w_1, w_2) \neq (0,0)\). Up to exchanging \(x\) and \(y\), we may assume \(w_1 \neq 0\) and set \(\gamma\defeq w_2/w_1\). For \(x>0\), factoring out \(x^{d/w_1}\) gives
		\[
			P_{xy}(x,y)
			=x^{d/w_1}\sum_{(m,n) \in \Lambda(P_{xy})} a_{mn} (yx^{-\gamma})^n
			\eqdef x^{d/w_1}Q(yx^{-\gamma}),
		\]
where \(Q \in \R[z]\setminus\set{0}\). Factoring \(Q(z)=uz^N\prod_{\ell=1}^L (z-c_\ell)^{k_\ell}\) over \(\C\) and substituting \(z=yx^{-\gamma}\) yields the desired form, where \(M \in \Z_{\ge 0}\) is forced by the polynomiality of \(P_{xy}\).
}

\section{Sobolev regularity of the global Radon transform} \label{appendix:radon}

In this appendix, we provide the proof of \Cref{cor:global_Radon} and demonstrate the failure of isotropic Sobolev regularity for the global Radon transform.

\Pf{Proof of \Cref{cor:global_Radon}}{
Let \(\widehat g(\lambda,y)\) denote the partial Fourier transform of \(g(s,y)\) in \(s\). Taking the Fourier transform in \(t\) gives
	\[
		\widehat{\cR g}(\lambda, x) = \cT^\lambda_\glo \paren[\big]{\widehat g(\lambda, \cdot)}(x).
	\]
Plancherel's theorem then yields
	\[
		\norm{\cR g}_{\dot H^{(\tau,0)}}^2
		= \int_\R |\lambda|^{2\tau}
			\norm[\big]{\cT^\lambda_\glo \paren[\big]{\widehat g(\lambda,\cdot)}}_{L_x^2}^2 \d{\lambda}.
	\]
Consequently, localizing \(\widehat g(\lambda, \cdot)\) in the \(\lambda\)-variable yields
	\begin{equation} \label{eq:Radon_fiber_norm}
		\norm{\cR}_{L^2(\R^2) \to \dot H^{(\tau,0)}(\R^2)} < + \infty
		\qquad \text{if and only if} \qquad
		\esssup_{\lambda \neq 0} |\lambda|^\tau \norm{\cT^\lambda_\glo}_{\op} < +\infty.
	\end{equation}

Suppose first that every Newton pair arising from the four quadrant pieces is balanced. By \Cref{cor:full_operator},
	\[
		\norm{\cT^\lambda_\glo}_{\op} \approx
		\begin{cases}
			|\lambda|^{-\frac{1}{2(d_1+1)}}, & \text{if } |\lambda| \ge 1, \\[1ex]
			|\lambda|^{-\frac{1}{2(d_2+1)}}, & \text{if } 0< |\lambda| < 1.
		\end{cases}
	\]
It follows that the essential supremum in \eqref{eq:Radon_fiber_norm} is finite if and only if
	\[
		\frac{1}{2(d_2+1)} \le \tau \le \frac{1}{2(d_1+1)}.
	\]
This proves (A).

Finally, suppose that some quadrant contains an unbalanced Newton pair. By \Cref{cor:full_operator}, \(\cT^\lambda_\glo\) is unbounded on \(L^2(\R)\) for every \(\lambda \neq 0\). Thus the essential supremum in \eqref{eq:Radon_fiber_norm} is infinite for every \(\tau \in \R\), proving (B).
}

\paragraph{Failure of regularity in the \(x\)-variable.}
For \(\tau, \sigma \in \R\), let \(\dot H^{(\tau,\sigma)}(\R^2)\) be the anisotropic homogeneous Sobolev space equipped with the seminorm
	\[
		\norm{F}_{\dot H^{(\tau,\sigma)}(\R^2)} 
		\defeq \norm[\big]{|D_t|^\tau |D_x|^\sigma F}_{L^2(\R^2)}.
	\]
The model \(P(x,y)=xy\) demonstrates that a bounded estimate with nonzero Sobolev order in the \(x\)-variable is impossible in general. We make this explicit by determining all pairs \((\tau,\sigma)\in \R^2\) for which \(\norm{\cR g}_{\dot H^{(\tau,\sigma)}(\R^2)} \lesssim \norm{g}_{L^2(\R^2)}\) holds. Continuing to write \(\widehat g(\lambda,y)\) for the partial Fourier transform of \(g(s,y)\) in \(s\), we take the Fourier transform in both \(t\) and \(x\) and obtain
	\[
		\widehat{\cR g}(\lambda, \xi) = c|\lambda|^{-1} \widehat g(\lambda, \xi/\lambda).
	\]
Consequently,
	\[
		\begin{alignedat}{2}
			\norm{\cR g}_{\dot H^{(\tau,\sigma)}}^2
			&\approx 
			\iint |\lambda|^{2\tau-2}|\xi|^{2\sigma}|\widehat g(\lambda, \xi/\lambda)|^2 \d{\xi}\d{\lambda} \\
			&=
			\iint |\lambda|^{2(\tau+\sigma)-1}|w|^{2\sigma} |\widehat g(\lambda,w)|^2 \d{w}\d{\lambda}.
		\end{alignedat}
	\]
If \(\sigma>0\), taking \(\widehat g(\lambda,w) = \chi(\lambda)\psi(w-N)\) and letting \(N \toinf\) shows that no such estimate holds. If \(\sigma<0\), the same conclusion follows by taking \(\widehat g(\lambda, w) = \chi(\lambda)\chi(Nw)\) and letting \(N\toinf\). Thus no bounded estimate holds for \(\sigma \neq 0\). 

When \(\sigma=0\), the remaining factor is \(|\lambda|^{2\tau-1}\). Uniform boundedness at both low and high frequencies therefore forces \(\tau=1/2\), in which case 
	\[
		\norm{\cR g}_{\dot H^{(1/2,0)}}^2 \approx \norm{g}_{L^2(\R^2)}^2.
	\]
Therefore, for \(P(x,y)=xy\), the global Radon transform \(\cR\) extends boundedly from \(L^2(\R^2)\) to \(\dot H^{(\tau,\sigma)}(\R^2)\) if and only if \((\tau,\sigma)=(1/2,0)\).


\printbibliography

\end{document}